\documentclass[11pt] {amsart}
\usepackage{amsthm}
\usepackage{amsxtra}
\usepackage{amssymb}
\usepackage{graphicx}

\newtheorem{theorem}{Theorem}
\newtheorem{prop}{Proposition}
\numberwithin{prop}{section}
\newtheorem{corol}{Corollary}
\numberwithin{corol}{section}

\newtheorem{lemma}{Lemma}
\numberwithin{lemma}{section}

{\theoremstyle{definition}
\newtheorem{defin}{Definition}
\numberwithin{defin}{section}
}
\numberwithin{figure}{section}

\renewcommand{\Re}{\mathop{\rm Re}\nolimits}
\renewcommand{\Im}{\mathop{\rm Im}\nolimits}

\newcommand{\ch}{\operatorname{ch}}
\newcommand{\WF}{\operatorname{WF}}
\newcommand{\ES}{\operatorname{ES}}
\newcommand{\supp}{\operatorname{supp}}

\newcommand{\CIc}{C^\infty_{\rm{c}}}
\newcommand{\RR}{\mathbb R}

\newcommand{\bl}{\begin{flushleft}}
\newcommand{\el}{\end{flushleft}}
\newcommand{\br}{\begin{flushright}}
\newcommand{\ert}{\end{flushright}}
\newcommand{\bc}{\begin{center}}
\newcommand{\ec}{\end{center}}
\newcommand{\mcal}[1]{\mathcal{#1}}

\newcommand{\recip}[1]{\frac{1}{#1}}

\newcommand{\complex}{\mathbb{C}}

\newcommand{\numList}{\begin{enumerate}}
\newcommand{\enumList}{\end{enumerate}}

\newcommand{\composed}{\text{\textopenbullet}}

\newcommand{\e}{\epsilon}

\newcommand{\cor}{\textbf{\begin{lemma}:}}

\newcommand{\re}{\mathbb{R}}

\newcommand{\nn}{\nonumber\\}
\newcommand{\la}{\langle}
\newcommand{\ra}{\rangle}
\newcommand{\negint}{{\int\negthickspace\negthickspace\negthickspace\negthinspace -}}
\newcommand{\cald}{Calder\'{o}n }

\numberwithin{equation}{section}

\usepackage{color}
\usepackage{amsfonts}
\usepackage{textcomp}
\usepackage{lmodern}

\title [Pseudospectra of Semiclassical Boundary Value Problems]{Pseudospectra of Semiclassical Boundary Value Problems}
\author[J. Galkowski]{Jeffrey Galkowski}
\address{Mathematics Department, University of California, Berkeley, 
CA 94720, USA}
\email{jeffrey.galkowski@math.berkeley.edu}

\begin{document}
\begin{abstract}
We consider operators $ - \Delta + X $, where $ X $ is a constant
vector field, in a bounded domain and show spectral
instability when the domain is expanded by scaling. More generally, we consider
semiclassical elliptic boundary value problems which exhibit spectral
instability for small values of the semiclassical parameter $h$, which
should be thought of as the reciprocal of the Peclet constant. This instability is due to the
presence of the boundary: just as in the case of  $ - \Delta + X $, 
some of our operators are normal when considered on $\re^d$. We
characterize the semiclassical pseudospectrum of such problems as well
as the areas of concentration of quasimodes. As an application, we prove
a result about exit times for  diffusion processes in bounded
domains.  We also demonstrate instability for a class of spectrally stable nonlinear
evolution problems that are associated to these elliptic operators.
\end{abstract}
\maketitle

\section{Introduction}
For many non-normal operators, the size of the resolvent is a measure of spectral
instability and is not connected to the distance to the spectrum. The sublevel sets of the norm of the
resolvent are referred to as the {\em pseudospectrum}.
The study of the pseudospectrum has been a topic of interest both in applied mathematics
(see \cite{Davi98},\cite{ET}, \cite{Hansen}, and numerous references given there) and
the theory of 
partial differential equations (see, for example
\cite{zworski01},\cite{De},\cite{KPS},\cite{Gall}, \cite{HiPS2}, \cite{HiPS1}). 

The problem of characterizing pseudospectra for semiclassical partial differential operators acting on Sobolev spaces on
$\re^d$ started with \cite{Davi98}. Dencker, Sj\"{o}strand, and
Zworski gave a more complete characterization for these pseudospectra in \cite{zworski04} by proving that, for operators with Weyl symbol $p$, if $(p-z)(x_0,\xi_0)=0$ and $i\{p,\bar{p}\}(x_0,\xi_0)<0$ then $z$ is in the semiclassical pseudospectrum of $p^w$. That is, for all $N>0$, there exists $C_N>0$ such that 
$$\left\|(p^w-z)^{-1}\right\|_{L^2\to L^2}\geq C_Nh^{-N}.$$
Moreover, they show that there exists a quasimode at $z$ in the following sense: there exists $u\in H_h^2$ (see \eqref{eqn:semiclassicalSobolev} for the definition of $H_h^2$) with $\|u\|_{L^2}=1$ and $\WF_h(u)=\{(x_0,\xi_0)\}$ such that 
$$(p^w-z)u=O_{L^2}(h^\infty).$$ In \cite{KPS}, Pravda-Starov extended the results of \cite{zworski04} and gave a slightly different notion of semiclassical pseudospectrum. 

In this paper, we examine the size of the resolvent for operators defined on bounded domains $\Omega\subset \re^d$ with $C^\infty$ boundary. Let \begin{equation}P=(hD)^2+i\la X,hD\ra\quad D_j:=\recip{i}\partial_j.\label{eqn:P}\end{equation}
where $X\in \re^d\setminus\{0\}.$
Here, $h$ can be thought of as the inverse of the Peclet constant. 
 We are interested in determining the semiclassical pseudospectrum of the Dirichlet operator $P$ on $\Omega$. That is, we wish to find $z\in \complex$ and $u\in H_h^2$ such that 
\begin{equation}
\label{eqn:generic}
\begin{cases}
P_zu:=(P-z)u=O_{L^2}(h^\infty)&x\in \Omega,\\
u|_{\partial\Omega}=0\,,\|u\|_{L^2}=1.
\end{cases}
\end{equation}

The collection of such $z$ will be denoted $\Lambda(P,\Omega)$ and the {\em{semiclassical pseudospectrum}} of $(P,\Omega)$ is $\overline{\Lambda(P,\Omega)}$. From this point forward, we will refer to $\overline{\Lambda(P,\Omega)}$ as the pseudospectrum. A solution to  (\ref{eqn:generic}) will be called a quasimode for $z$. We restrict our attention to the case where $X$ is constant so that there are no quasimodes given by the results of \cite{zworski04} and, moreover, the operator is normal when acting on $L^2(\re^d)$.

We characterize $\overline{\Lambda(P,\Omega)}$ for such boundary value problems as well as the semiclassical essential support of quasimodes. Here the essential support is defined as 

\begin{figure}[htbp]

\centering
\def\svgwidth{3in}
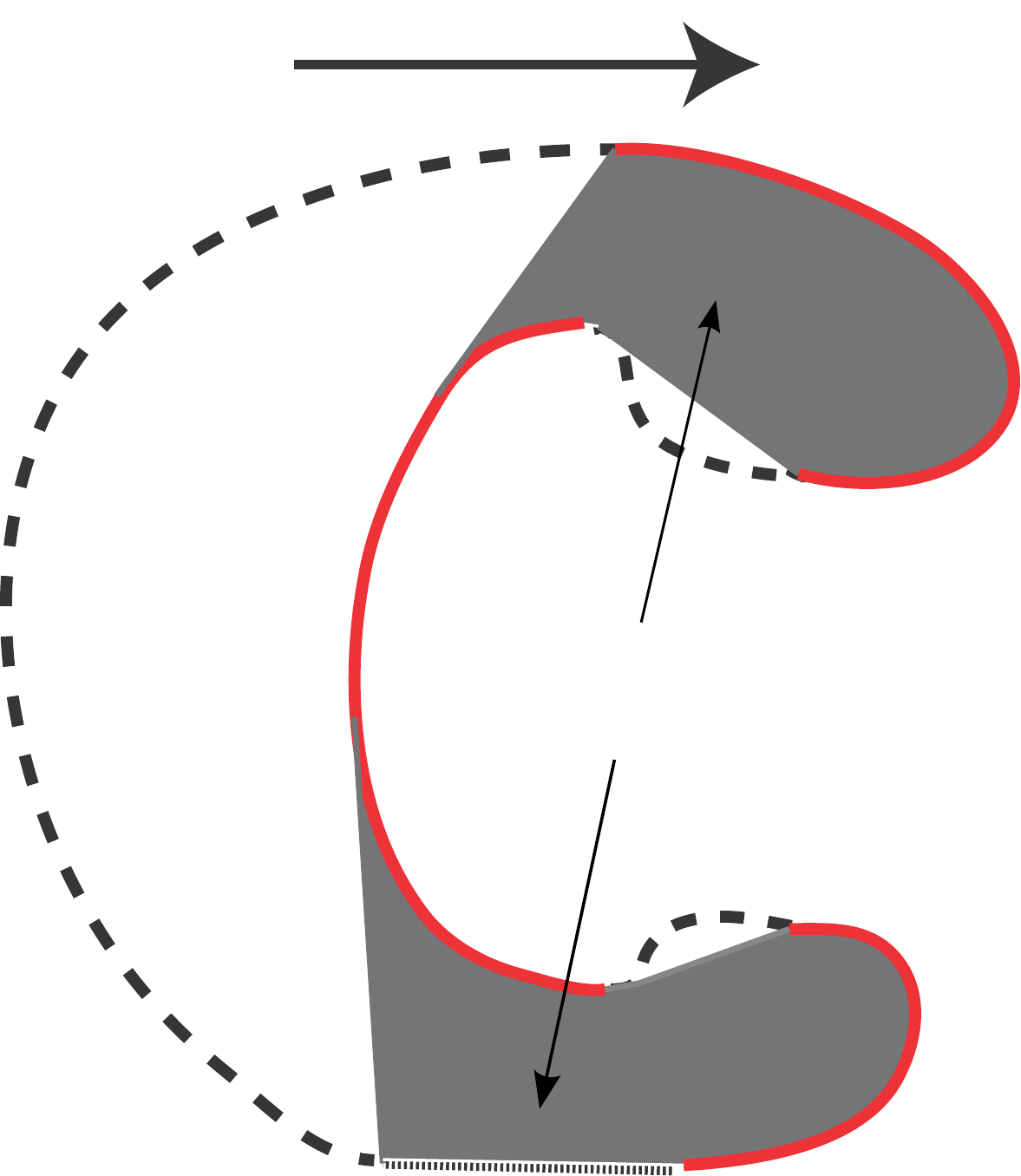
\caption{\label{f:Omega}
The figure shows an example domain $\Omega\subset \re^2$. $\partial\Omega_+$ is shown in the solid red line, $\partial\Omega_0$ in the dotted black line, and $\partial\Omega_-$ in the dashed black line. The region $\overline{\ch_{\overline{\Omega}}(\Gamma_+)}$ is shaded in. The $\psi$ shown here is the locally convex function used to prove Lemma \ref{lem:convexRestriction}}
\end{figure}

\begin{defin} The essential support of a family of $ h$-dependent
  functions $u = u ( h ) $ is given by
$$\ES_h(u):=\bigcap_{U\in \mcal{A}} U \,, \quad \mcal{A}:=\{U\subset\overline{\Omega}:\text{ if }\chi\in C^\infty(\overline{\Omega}),\text{ }\chi\equiv 1 \text{ on }U,\text{ then }(1-\chi)u=O_{L^2}(h^\infty)\}.$$
\end{defin}
We will need the following analogue of convexity (similar to that used for planar domains in \cite{relConv}).
First, define 
\begin{equation}
\label{eqn:lineSeg}
L_{x,y}:=\{tx+(1-t)y:t\in[0,1]\},
\end{equation}
the line segment between $x$ and $y$. 
Then,
\begin{defin} A set $A\subset B$ is relatively convex in $B$ if for all $x,y\in A$, $L_{x,y}\subset B$ implies $L_{x,y}\subset A$. 
\end{defin}
We also need an analogue of the convex hull in this setting
\begin{defin}
\label{def:ConvSub}
For $A\subset B$, we define the convex hull of $A$ relative to $B$ by 
$$\ch_B(A)=\bigcap_{C\in\mcal{A}}C \,, \quad \mcal{A}:=\{C:A\subset C\text{ and }C \text{ is relatively convex in }B\}.$$
(See Figure \ref{f:Omega} for an example.)

\noindent If $A\nsubseteq B$, then define 
$$\ch_B(A)=\ch_B(A\cap B).$$
\end{defin}

\noindent{\bf Remark:} In the case that $B$ is convex, these definitions coincide with the usual notions of convexity.

Let $\nu$ be the outward pointing unit normal to $\partial \Omega$. We define subsets of $\partial \Omega$ similar to those in \cite{SjoUhl},
\begin{equation}
\partial{\Omega_-}=\{x\in \partial \Omega:\la X, \nu\ra< 0\},\quad\label{eqn:OmegaBound}
\partial{\Omega_+}=\{x\in \partial \Omega:\la X, \nu\ra> 0\},\end{equation}
\begin{equation}\nonumber
\partial\Omega_0=\{x\in \partial \Omega: \la X,\nu\ra =0\}, \quad \Gamma_+=\partial\Omega_+\cup\partial\Omega_0.
\end{equation}

\noindent{\bf Remark: }We will refer to $\partial\Omega_+$, $\partial\Omega_0$, and $\partial\Omega_-$ as the illuminated, glancing, and shadow sides of the boundary, respectively. Figure \ref{f:Omega} shows examples of these subsets in a two dimensional domain.

With these definitions in place, we can now state our main theorem:
\begin{theorem}
\label{thm:PseudoSpec}
Let $P$ be as in (\ref{eqn:P}), and $\Omega\subset\re^d$ be a domain with $C^\infty$ boundary. Then, 
\begin{enumerate}
\item $\overline{\Lambda(P,\Omega)}=\{z\in \complex: \Re z\geq (\Im z)^2|X|^{-2}\}.$
 (Here $|\cdot|$ is the Euclidean norm.) 
\item For all quasimodes $u$,
$$\ES_h(u)\subset\overline{ \ch_{\overline{\Omega}}(\Gamma_+)}\cap \partial\Omega \quad \text{ and }\quad \ES_h(u)\cap\overline{\partial\Omega_+}\neq\emptyset.$$ 
\item For each point $x_0\in\partial\Omega_+$, there exists
$$W_{x_0}\subset\{z\in \complex:\Re z> (\Im z)^2|X|^{-2}\}$$
 such that $W_{x_0}$ is open and dense in $\overline{\Lambda(P,\Omega)}$ and for each $z\in W_{x_0}$, there is a quasimode $u$ for $z$ with $\ES_h(u)=x_0$. Moreover, if $\partial\Omega$ is real analytic near $x_0$, then there exists $c>0$ such that these quasimodes can be constructed with $P_zu=O_{L^2}(e^{-c/h})$. 
 \item Let $x_0\in \partial\Omega_-$. Suppose that $\partial\Omega$ is strictly convex or strictly concave at $x_0$. Then, for any quasimode $u$, $x_0\notin \ES_h(u).$
\item  If $\Omega\subset \re^2$, and $u$ is a quasimode, then $\ES_h(u)\subset \Gamma_+.$
\end{enumerate}
\end{theorem}

\noindent {\bf Remark:} If $\Omega$ is convex then Theorem \ref{thm:PseudoSpec} gives that $\ES_h(u)\subset\Gamma_+.$

\begin{figure}[htbp]
\includegraphics[width=6in]{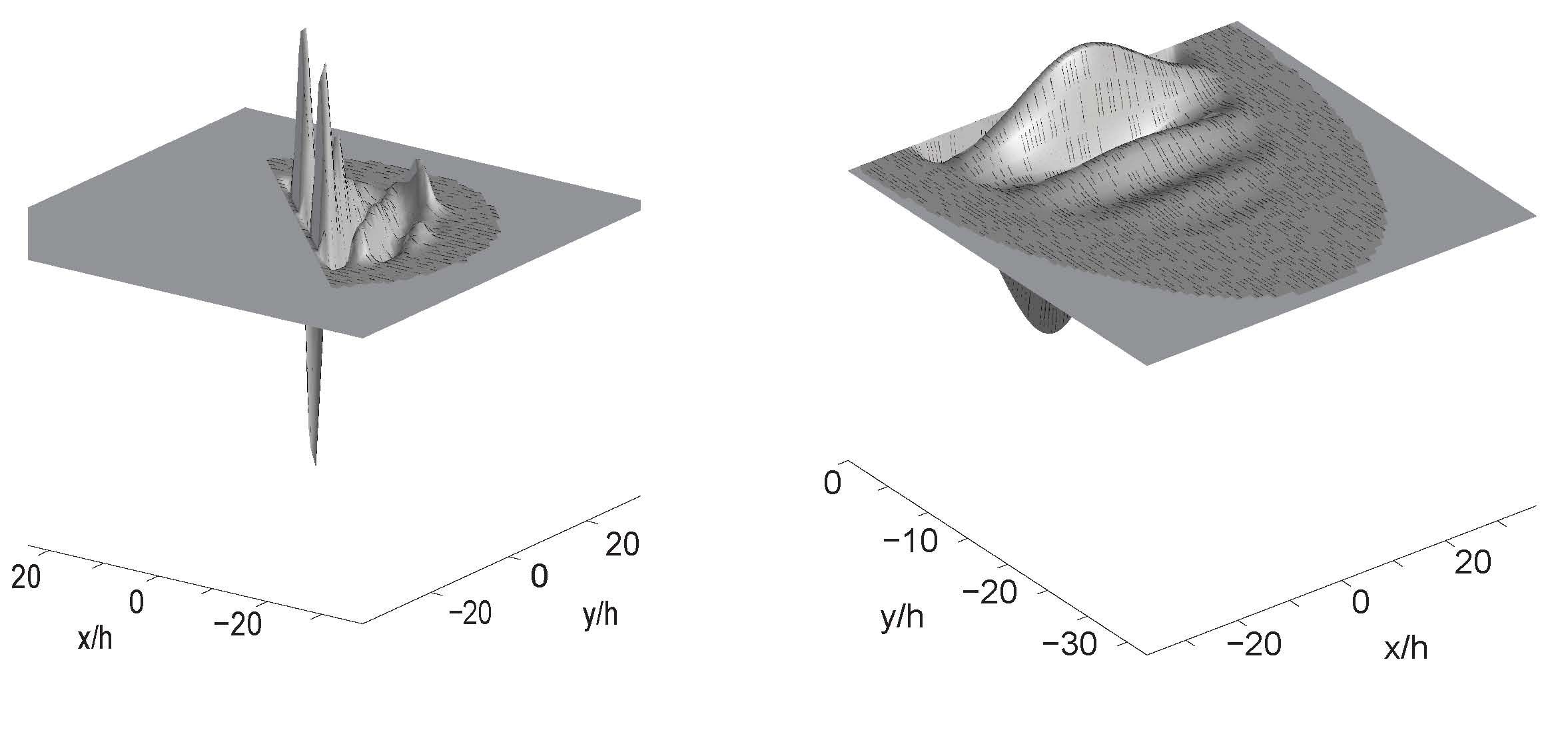}
\begin{center}
\caption{\label{f:quasimode}
The figure shows an example of the imaginary part of a quasimode constructed in Proposition \ref{prop:constructQuasi}. On the left, the boundary forms an angle of $\pi/4$ with $X=\partial_{x_1}$, and on the right the boundary is normal to $\partial_{x_1}$. In both cases, $ z=1+i/2$.}
\end{center}
\end{figure}

\begin{figure}[htbp]
\includegraphics[width=3.5in]{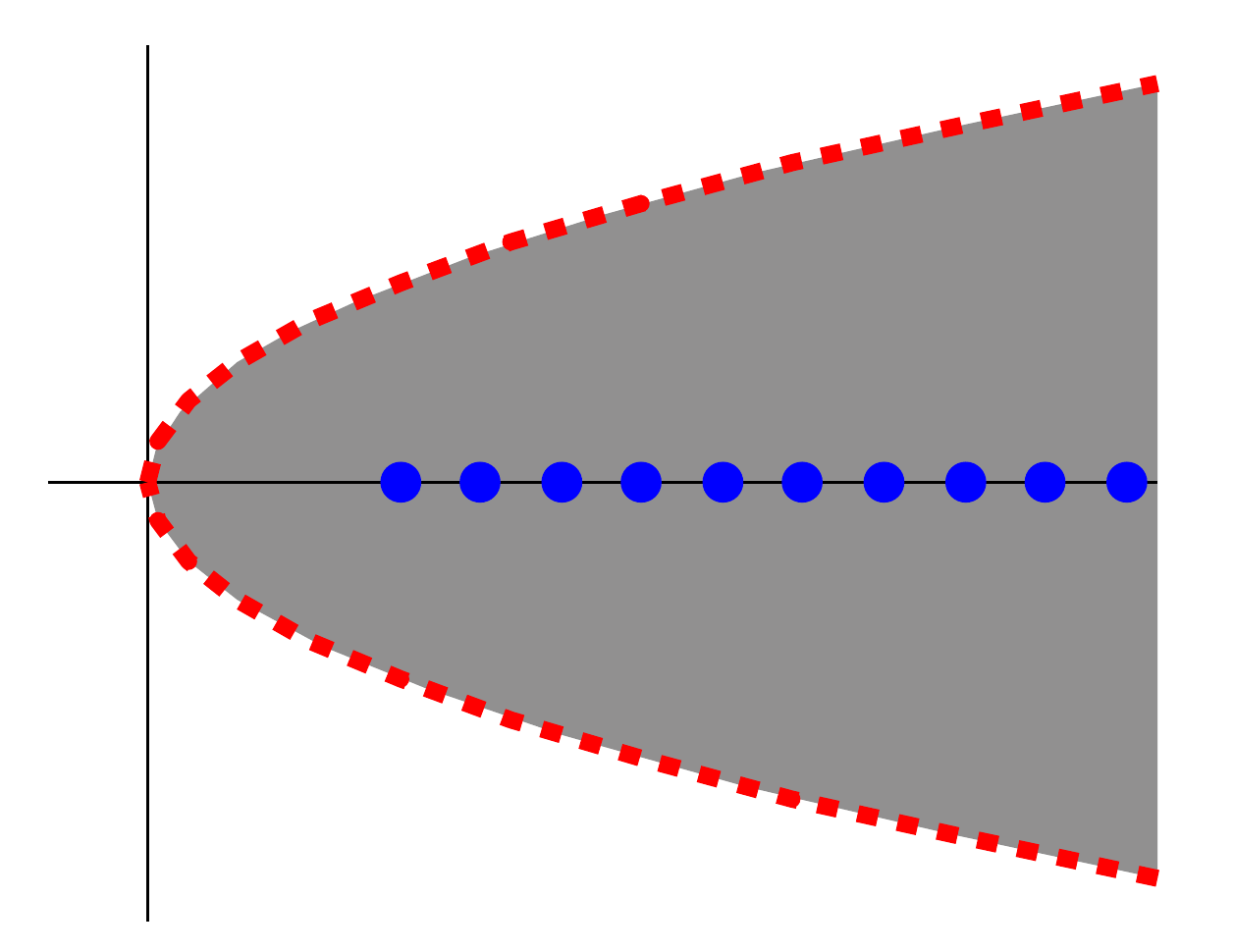}
\begin{center}
\caption{\label{f:Pseud}
The figure shows the pseudospectrum and spectrum of $(P,\Omega)$ for $|X|=1$. The pseudospectrum is the shaded region, the spectrum is shown as blue circles, and the curve $\Re z=(\Im z)^2$ is shown in dashed red. The spectrum of $(P,\Omega)$ is discrete and real since $P$ is an elliptic second order partial differential operator. Moreover, in the case that $X$ is a constant, $P$ can be conjugated to a self-adjoint elliptic operator using a non-unitary operator and hence has real spectrum.}
\end{center}
\end{figure}

When $X\neq 0$ is constant, conjugating $P$ by $e^{-\la X,x\ra/2h}$ shows that the spectrum of $(P,\Omega)$ is discrete and contained in $\{z\in \complex: \Re z\geq c>0,\text{ }\Im z=0\}$. Thus, Theorem \ref{thm:PseudoSpec} shows that the pseudospectrum of $(P,\Omega)$ is far from its spectrum and hence that the size of the resolvent is unstable in the semiclassical limit. (Figure \ref{f:Pseud} shows the spectrum and pseudospectrum of $(P,\Omega)$ in an example.)

For a large class of nonlinear evolution equations this type of behavior has been proposed as an explanation of instability 
for spectrally stable problems. Celebrated examples include the 
plane Couette flow, plane Poiseuille flow and plane flow -- 
see Trefethen-Embree \cite[Chapter 20]{ET} for discussion and
references. Motivated by this, we consider the mathematical question 
of evolution involving a small parameter $ h $ (in fluid
dynamics problem we can think of $ h $ as the reciprocal of
the Reynolds number) in which the linearized operator has spectrum lying in $ \Re z <  - \gamma_0 < 0 $, uniformly 
in  $h $,  yet the solutions of the nonlinear equation blow up in short time for data of size $ O ( \exp ( - c/h )) $.   

Let $p>1$. We examine the behavior of the following nonlinear evolution problem
\begin{equation}\label{eqn:evolutionEq}
\begin{cases}
(h\partial_t+(P-\mu))u-u^p=0,& t\geq 0,\quad x\in \Omega\subset \re^d,\\
u|_{\partial\Omega}=0&u(x,0)=u_0(x),\end{cases}
\end{equation}
and interpret it in terms of the pseudospectral region of $P-\mu$.

We have the following analog of what is shown in \cite{Ga} and \cite{Sand03}
\begin{theorem}\label{thm:evolution}
Fix $\mu>0$. Then, for 
\[ 0 < h < h_0 \,, \]
where $ h_0 $ is small enough, and each $\delta>0$ , there exists
\[ u_0 \in \CIc ( \RR^n ) \,, \ \ u_0 \geq 0 \,, \ \ \| u_0 \|_{ C^k} \leq
\exp\left(-\recip{C_kh}\right) \,, \  \ k = 0 , 1 , \dots \,, \]
such that the solution to \eqref{eqn:evolutionEq} with $ u ( x, 0 ) =  u_0 (
x ) $, satisfies 
\[   \| u ( x, t ) \|_{L^\infty } \longrightarrow \infty , \ \ t
\longrightarrow T \,, \]
where $   T \leq \delta \, . $
\end{theorem} 

As an application of Theorem \ref{thm:PseudoSpec}, we consider diffusion processes on bounded domains. Specifically, we examine hitting times
$$\tau_X=\inf\{t\geq 0:X_t\in \partial \Omega\}$$
 for processes of the form 
$$dX_t=b(X_t)+\sqrt{2h}dB_t\quad X_0=x_0(h)$$
where $B_t$ is standard Brownian motion in $d$ dimensions, $x_0(h)\to x_0\in \partial\Omega_+$ (defined for the vector field $-b$), and $b\in C^\infty(\re^d;\re^d)$. We show that, for $\partial\Omega_+$ analytic near $x_0$, and all $N\geq1$, the log moment generating function of $\tau_x$ does not decay as $h\to 0$ for $|x_0(h)-x_0|\approx Ch^N$. Moreover, letting 
$$L:=(hD)^2+i\la -b,hD\ra$$
 and $0<c<\lambda_1(L)$ be the principal eigenvalue of $L$, we have for each $0<\lambda<\lambda_1(L)$ (where $\lambda$ is $h$ independent), each $\e>0$, and $x_0(h)$ as above, that there exists $\delta>0$ such that for all $\alpha>1$, there exists a function $s(h)>\delta-h^{1-\e}$ and $c_\alpha>0$ a constant depending only on $\alpha$ such that for $h$ small enough,
$$\min(c_\alpha e^{-\alpha (s(h)-\delta)/h},1) \leq P\left(\tau_X\geq \frac{s(h)}{\lambda}\right)\leq P\left(\tau_X\geq \frac{\delta-h^{1-\e}}{\lambda}\right).$$

\noindent {\bf Remark:} See the remarks after Proposition \ref{prop:exitTime} for an interpretation of this inequality.

\noindent
{\sc Acknowledgements.} The author would like to thank Maciej Zworski for suggesting the problem and for valuable discussion. Thanks also to Fraydoun Rezakhanlou for advice on the application to diffusion processes and to the anonymous referee for many helpful comments. The author is grateful to the National Science Foundation for partial support under grant DMS-0654436 and under the National Science Foundation Graduate Research Fellowship Grant No. DGE 1106400.
\section{Outline of the Proof}
In this section, we explain the ideas of the proof of Theorem \ref{thm:PseudoSpec}. We also describe the structure of the paper. 

Our starting point is to prove that if $|p_z(x,\xi)|>C\la \xi \ra^2$, for all $x\in \overline{\Omega}$, then \eqref{eqn:generic} has an inverse that is bounded independently of $h$ on semiclassical Sobolev spaces. We do this via a construction of \cald projectors adapted from \cite[Chapter 20]{HOV3}. It follows from the existence of such inverses that 
$$\Lambda(P,\Omega)\subset \{\Re z\geq (\Im z)^2|X|^{-2}\}.$$

Next, we show that $$\overline{\Lambda(P,\Omega)}=\{\Re z\geq(\Im z)^2|X|^{-2}\}.$$ In particular, we construct quasimodes near points $x_0\in \partial \Omega_+$. To do this, we use a WKB method adapted to Dirichlet boundary value problems. Motivated by the fact that, in one dimension, eigenfunctions of the Dirichlet realization of $P$ are of the form $e^{cx/h}\sin(x/h)$, we look for solutions of the form 
$$a(x)e^{i\varphi_1(x)/h}-b(x)e^{i\varphi_2(x)/h}$$
and derive formulae for WKB expansions of $a$ and $b$. In order to complete this construction, we have to solve a complex eikonal equation for the $\varphi_i$'s. This is done by finding the Taylor expansions $\varphi_i$ at $x_0$. We proceed similarly for $a$ and $b$. Figure \ref{f:quasimode} shows examples of quasimodes constructed using this method.

Our last task is to characterize the essential support of quasimodes $u$. The main idea is to prove a Carleman type estimate for solutions to \eqref{eqn:generic}. This estimate gives us control of solutions outside relatively convex sets containing $\Gamma_+$. Hence any quasimode is essentially supported inside such a convex set. The next ingredient in the proof is a result adapted from \cite{FIO2} on propagation of semiclassical wavefront sets for solutions of \eqref{eqn:generic}. We show that the wavefront set of a quasimode is invariant under the leaves generated by $H_{\Im p}=\la X,\partial_x\ra$ and $H_{\Re p}.$ We then show that there exist convex sets containing $\Gamma_+$ which do not extremize $\la X,x\ra$ inside $\Omega$. Finally, we combine this with the propagation results to show that $$\ES_h(u)\subset \overline{\ch_{\overline{\Omega}}(\Gamma_+)}\cap \partial\Omega.$$ 

The paper is organized as follows. In section \ref{sec:notation} we introduce various semiclassical notations. Then, in section \ref{sec:cald} we prove results on \cald projectors adapted from \cite[Chapter 20]{HOV3}. Section \ref{sec:QuasiConstruct} contains the construction of quasimodes via a boundary WKB method. In section \ref{sec:propagation}, we adapt results of Duistermaat and H\"{o}rmander in \cite[Chapter 7]{FIO2} on propagation of wavefront sets to the semiclassical setting. In section \ref{sec:Carleman}, we prove a Carleman type estimate that will be used in section \ref{sec:QuasimodesOnIlluminated} to derive restrictions on the essential support of quasimodes.  Section \ref{sec:Instabilty} contains the proof of Theorem \ref{thm:evolution}. Finally, section \ref{sec:meanHit} applies some of the results of Theorem \ref{thm:PseudoSpec} to exit times for diffusion processes.

\section{Semiclassical Preliminaries and Notation}
\label{sec:notation}
The $ O(\cdot)$ and $o(\cdot)$ notations are used in the present paper in the following ways:
we write $u=O_{\mathcal X}(F(s))$ if the norm of the function,
or the operator, $u$ in the functional space $\mathcal X$ is bounded
by the expression $F$ times a constant independent of $s$. We write $u=o_{\mcal{X}}(F(s))$ if the norm of the function or operator, $u$ in the functional space $\mcal{X}$ has 
$$\lim_{s\to s_0}\frac{\|u\|_{\mcal{X}}}{F(s)}=0$$
where $s_0$ is the relevant limit. If no space $\mcal{X}$ is specified, then this is understood to be pointwise. 

The \emph{Kohn-Nirenberg symbols} for $m\in \mathbb{R}$ as in \cite[Section 9.3]{EZB} by
$$S^m:=\{a(x,\xi)\in C^\infty(\re^{2d}): \left|\partial_x^\alpha\partial_\xi^\beta a\right|\leq C_{\alpha\beta}\la \xi\ra^{m-|\beta|}\},\quad \la \xi \ra=\left(1+|\xi|^2\right)^{1/2}$$
and denote by $\Psi^m$, the semiclassical pseudodifferential operators of order $m$, given by 
$$\Psi^m:=\{a(x,hD)|a\in S^m\}.$$

\noindent{\bf Remark:} We will sometimes write $S=S^0$ and $h^kS^m$ denotes the class of symbols in $S^m$ whose seminorms are $O(h^k)$.

Throughout this paper, we will use the standard quantization for pseudodifferential operators on $\re^d$, 
$$a(x,hD)u=(2\pi h)^{-d}\int\!\!\!\!\int a(x,\xi)e^{\frac{i\la x-y,\xi\ra}{h}}u(y)dyd\xi,$$
unless otherwise stated. In those cases that we use the Weyl quantization,
$$a^w(x,hD)u=(2\pi h)^{-d}\int\!\!\!\!\int a\left(\frac{x+y}{2},\xi\right)e^{\frac{i\la x-y,\xi\ra}{h}}u(y)dyd\xi,$$
the operator will be denoted $p^w(x,hD)$ where $p$ is the symbol of the operator. Using semiclassical pseudodifferential operators, we can now define the \emph{semiclassical Sobolev spaces} $H_h^s:=\la hD\ra ^{-s}L^2$ with norm 
\begin{equation}
\label{eqn:semiclassicalSobolev}
\|u\|_{H_h^s}=\|\la hD\ra^s u\|_{L^2}.
\end{equation}
For more details on the calculus of pseudodifferential operators see, for example, \cite[Chapter 4]{EZB}.

We briefly recall the definition of pseudodifferential operators on a compact manifold $M$. We say that an operator $B:\mcal{S}(M)^n\to \mcal{S}'(M)^k$ is a pseudodifferential operator, denoted $B\in \Psi^m(M,\mathbb{C}^k\otimes\mathbb{C}^n)$, if, 
\begin{enumerate}
\item letting $U_\gamma$ be a coordinate patch, and $\varphi$,$\psi\in C_0^\infty(U_\gamma)$, the kernel of $\varphi (\gamma^*)B(\gamma^{-1})^*\psi$ can be written as an $n\times k$ matrix $b_{ij}$ where $b_{ij}\in S^m$.
\item if $\varphi,\psi\in C_0^\infty(M)$ with supp $\psi\cap$ supp $\varphi=\emptyset$ then for all $N$
$$\varphi B\psi=O(h^\infty):H_h^{-N}(M;\complex^n)\to H_h^N(M;\complex^k).$$
\end{enumerate}
See \cite[Chapter 14]{EZB} for a detailed account of pseudodifferential operators on manifolds.

Finally, we need a notion of microlocalization for semiclassical functions. We call $u$ \emph{tempered} if for some $m,N>0$, 
$$u\in H_h^m,\quad\|u\|_{H_h^m}\leq Ch^{-N}.$$
For a tempered function, $u$, we define the \emph{semiclassical wavefront set of $u$}, $\WF_h(u)$, by $(x_0,\xi_0)\notin \WF_h(u)$ if there exists $a\in S^0$ with $|a(x_0,\xi_0)|>\gamma>0$ such that 
$$\|a^w(x,hD)u(h)\|_{H_h^m}=O(h^\infty).$$
(For more details on the semiclassical wavefront set see \cite[Section 8.4]{EZB}.)\\

\vbox{\noindent{\bf Other Notation: }
\begin{itemize}
\item Throughout the paper, we will denote the outward unit normal to $\partial \Omega$ at a point $x_0$, by $\nu(x_0)$. 
\item We will identify $T^*\re^d$ with $T\re^d$ using the Euclidean metric, denote by $|\cdot |$ the induced Euclidean norm on $T \re^d$, and $\la \cdot,\cdot \ra$ the inner product.
\item We will denote by $e_i$, the unit vector in the $x_i$ direction and $\nu_i=\la e_i,\nu\ra$. 
\item $\Omega^o$ will denote the interior of $\Omega$ and $\overline{\Omega}$, its closure.
\end{itemize}
}

\section{\cald  Projectors for Elliptic Symbols}
\label{sec:cald}
Our goal is to find an inverse, uniformly bounded in $h$, for the following elliptic boundary value problem
\begin{equation} \begin{cases}
\label{eqn:BVP}
Pu=f&\text{ in } \Omega,\\
B_ju=g_j,\text{ }j=1,...,J&\text{ on } \partial\Omega,
\end{cases}
\end{equation}
where $P=p(x,hD)$ is a differential operator of order $m$ with $p\in S^m$ elliptic (in the semiclassical sense, i.e. $|p|\geq c\la \xi\ra^m$), and $B_j$ are differential operators on the boundary of order $m_j$ with symbols $b_j\in S^{m_j}$.

\noindent{\bf Remark:}  In our applications, $B_1$ is the identity and $J=1$.

We define classical ellipticity for a boundary value problem as in \cite[Definition 20.1.1]{HOV3},
\begin{defin}
The boundary value problem \eqref{eqn:BVP} is called classically elliptic if 
\begin{enumerate}
\item For all $x\in \overline{\Omega}$, $|p(x,\xi)|\geq C\la \xi \ra ^m$ for $|\xi|\geq C_1$.
\item The boundary conditions are elliptic in the sense that for every $x\in \partial \Omega$ and $\xi \in T^*_x(\Omega)$ not proportional to the interior conormal $n_x$ of $X$, the map 
$$M_{x,\xi}^+\ni u\mapsto (b_1(x,\xi+D_tn_x)u(0),...,b_J(x,\xi+D_t n_x)u(0))\in \complex^J$$
is bijective, if $M_{x,\xi}^+$ is the set of all $u\in C^\infty(\re)$ such that $p(x,\xi+D_tn_x)u(t)=0$ for $t\in \re$ and such that $u$ is bounded on $\re_+$.
\end{enumerate}
\end{defin} 

\noindent We follow H\"{o}rmander's construction \cite[Chapter 20]{HOV3} of a \cald projector for classically elliptic boundary value problems to prove the following Proposition:

\begin{prop}
\label{prop:Calderon}
Let $P$ be as in \eqref{eqn:BVP}, $p(x,\xi)$ have $\inf_{x\in \overline{\Omega}}|p(x,\xi)|\geq c\la \xi \ra^m$, and \eqref{eqn:BVP} be classically elliptic. Then, for $h$ small enough and $s\geq m$, the system \eqref{eqn:BVP} has an inverse
 $$P_\Omega^{-1}:H_h^{s-m}(\Omega) \oplus_{j=0}^{m-1}H_h^{s-m_j-\recip{2}}(\partial\Omega)\to H_h^{s}(\Omega)$$ 
 with $\|P_{\Omega}^{-1}\|\leq C$ uniformly in $h$. 
\end{prop}

\subsection{Pseudospectra Lie Inside the Numerical Range}
Observe that Proposition \ref{prop:Calderon} gives that pseudospectra for elliptic boundary value problems must lie inside the numerical range of $p(x,\xi)$. In the special case of $P$ as in \eqref{eqn:P}, we have that $P=p(x,hD)$ where $p(x,\xi)=|\xi|^2+i\langle X,\xi\ra$. By Proposition \ref{prop:Calderon}, we have that, if $P_z$ is strongly elliptic, i.e. $|p_z(x,\xi)|\geq c\langle \xi\rangle^2$, then no quasimodes for $z$ exist.

Using this, observe that $p_z(x_0,\xi_0)=0$ implies
$$
\langle X,\xi_0\rangle =\Im z,\quad |\xi_0|^2-\Re z=0.
$$
Hence, identifying $T_{x_0}^*\re^d$ with $T_{x_0}\re^d$ using the Euclidean metric, there exists $w\in T_{x_0}^*\re^d$ with $\la X, w\ra=0$ such that 
\begin{equation}
\label{eqn:xiCond}
\xi_0=\Im z X|X|^{-2}+w,\quad
|w|^2=\Re z- |X|^{-2}\left(\Im z\right)^2.
\end{equation}
This, together with Proposition \ref{prop:Calderon}, implies that 
\begin{equation}
\label{eqn:ellipticRestriction}
 \Lambda(P,\Omega)\subset \{z\in \complex: \Re z \geq |X|^{-2}\left(\Im z\right)^2\}.
\end{equation}
\subsection{Proof of Proposition \ref{prop:Calderon}}

We will follow H\"{o}rmander's proof from \cite[Chapter 20]{HOV3} almost exactly. We present the proof in detail to provide a reference for \cald projectors in the semiclassical setting. Note that, unlike for operators that are only classically elliptic (in which case projectors yield a $C^\infty$ parametrix), in the semiclassically elliptic setting, the construction yields an inverse for the boundary value problem.

First observe that in \eqref{eqn:BVP} we can assume without loss of generality that the order of $B_j$ transversal to $\partial\Omega$ is less than $m$. To see this, let $U$ be a local coordinate patch in which $\Omega$ is defined by $x_1\geq 0$. Then, $Pu=f$ has the form 
$$\sum P_\alpha (x)(hD)^\alpha u =f$$
where $P_\alpha$ are $C^\infty$ functions of $x$. Then, observe that since $P$ is elliptic, the coefficient of $D_1^m$ is nonzero and we can write

$$(hD_1)^mu=P_{me_1}^{-1}(f-\sum_{\alpha_n<m}P_\alpha(x)(hD)^\alpha u),\quad e_1=(1,0,...,0).$$

\noindent Hence, if the transversal order of $B_j$ is greater than $m-1$, we can replace $(hD_1)^m$ by this expression. Using a partition of unity on the boundary to combine these local constructions, we obtain that 
$$B_j=B_j^r+C_jP$$
where $B_j^r$ has transversal order $<m$ and $C_j$ is a boundary differential operator. Then, \eqref{eqn:BVP} is equivalent to 
$$\begin{cases}
Pu=f&\text{ in }\Omega\\
B_j^ru=g_j-C_jf,\text{ } j=1,...,J&\text{ in }\partial \Omega
\end{cases}.
$$
 
Now, extend $P$ to a neighborhood, $\hat{\Omega}$ of $\overline{\Omega}$ so that $P$ is strongly elliptic on $\hat{\Omega}$. Then, define $T=P^{-1}$ where $T$ exists and is a pseudodifferential operator since $p$ is semiclassically elliptic i.e. $T$ is given by
$$T=\left[\recip{p}(x,hD)p(x,hD)\right]^{-1}\recip{p}(x,hD)\in \Psi^{-m}.$$ 
(The fact that $T$ is pseudodifferential follows from Beals' Theorem, (see for example \cite[Section 9.3.4 and remark after Theorem 8.3]{EZB}))
We will construct the \cald projector locally and hence reduce to the case where $\partial \Omega=\{x_1=0\}$ by a change to semigeodesic coordinates for $\partial\Omega$ and an application of a partition of unity. 

For $u \in C^\infty(\overline{\Omega})$, define 
$$\gamma u:=(u,(hD_1)^1u,...,(hD_1)^{m-1})|_{\partial\Omega}\in C^\infty(\partial \Omega).$$
In $\partial \Omega\times [0,1)$, we have 
$$P=\sum_{j=0}^mP_j(x,hD_2,...,hD_n) (hD_1)^j$$
where $P_j$ are semiclassical differential operators of order $m-j$ in $\partial \Omega$ depending on the parameter $x_1$. We denote the principal symbol of $P_j$ by $\sigma(P_j)=p_j$. Next, let $u^0$ denote extension by $0$ off of $\overline{\Omega}$. We have 
$$Pu^0=(Pu)^0+P^c\gamma u$$
where for $U=(U_0,...,U_{m-1})\in C^\infty (\partial \Omega)$, 
$$P^cU=i^{-1}\sum_{j<m}hP_{j+1}\sum_{k\leq j}U_{j-k}\otimes (hD_1)^k\delta,$$
where $\delta$ is the Dirac mass at $x_1=0$.
Then, 
\begin{equation}
\label{eqn:Tinv}
u^0=T[(Pu)^0+P^c\gamma u].\end{equation}

Next, define $Q$, the \cald projector, for $U\in C^\infty (\partial \Omega, \complex ^m)$, by $QU:=\gamma TP^cU$. Then, for $k=0,..., m-1$,
$$(QU)_k=\sum_0^{m-1}Q_{kl}U_l$$
with 
$$Q_{kl}U_l=\left.\left[\sum_{j=0}^{m-1-l}hi^{-1}(hD_1)^kTP_{j+l+1}U_l\otimes (hD_1)^j\delta\right]\right|_{\partial \Omega}.$$
(Note that the boundary values are taken from $\Omega^o$.) Therefore, we have that $Q_{kl}$ are pseudodifferential operators in $\partial \Omega$ of order $k-l$ with principal symbols
\begin{equation}
\label{eqn:Qsymb}q_{kl}(x',\xi ')=(2\pi i)^{-1}\int^+\sum_{j=0}^{m-l-1}\xi_1^{k+j}p(x',0,\xi)^{-1}p_{j+l+1}(0,x',\xi ')d\xi_1,\end{equation}
where the $\int^+$ denotes the sum of residues for $\Im \xi_1>0$.
\begin{lemma}
Q is a projection on the space of Cauchy data in the sense that $Q^2-Q=0$. If we identify solutions of the ordinary differential equation $p(0,x',hD_1,\xi ')v=0$ with the Cauchy data $(v(0),...,(hD_1)^{m-1}v(0))$ then, $q(x',\xi ')$ is for $\xi '\neq 0$ identified with the projection on the subspace $M^+$ of solutions exponentially decreasing on $\re_+$ along the subspace $M^-$ of solutions exponentially decreasing in $\re_-$. 
\end{lemma}
\begin{proof}
Let $u=TP^c U$,  by (\ref{eqn:Tinv}),
$$Pu=PTP^cU=P^cU=0\text{ in }\Omega^o.$$
Also, by (\ref{eqn:Tinv}),
\begin{align}
QU=\gamma (TP^c U)^0&=\gamma T((PTP^c U)^0+P^c\gamma u)=\gamma T(P^c U)^0+\gamma TP^c\gamma TP^cU = Q^2U\nonumber 
\end{align}
since $P^c U=0$ in $\Omega ^o$. 
Hence, $Q^2-Q=0$ and thus is a projection.

To see the second part of the claim, let $U$ as above. Then, the inverse semiclassical Fourier transform
$$v(x_1)=(2\pi h i)^{-1}h\int p(0,x',\xi)^{-1}\sum_{j+l<m}p_{j+l+1}(0,x',\xi ')\xi_1^jU_le^{ix_1\xi_1/h}d\xi_1$$
is in $\mcal{S}'$ and satisfies
$$p(0,x',\xi ',hD_1)v=hi^{-1}\sum_{j+l<m}p_{j+l+1}(x',0,\xi')U_l(hD_1)^j\delta$$ 
and hence $v$ coincides for $x_1>0$ with an element $v^+\in M^+$ and for $x_1<0$ with an element $v^-\in M^-$. For $x_1=0$, we have the jump condition 
$$(hD_1)^k(v^+-v^-)=U_k.$$

Then, (\ref{eqn:Qsymb}) gives 
$$q(x',\xi ')U=(v^+(0),...,(hD_1)^{m-1}v^+(0))$$
and $qU=U$ implies $v^-=0$, i.e. U Is the Cauchy data of a solution in $M^+$. Also, if $qU=0$, $U$ is the Cauchy data for $v\in M^-$ and we have proven the claim.
\end{proof}

Now that we have $Q$ defined locally, we can extend it to a global pseudodifferential operator on $\partial \Omega$ by taking a locally finite partition of unity, $\chi_j$ subordinate to $V_j$, the semigeodesic coordinate patches, and letting 
$$Q=\sum_j\chi_jQ_j.$$ 

To complete the proof, we need the following lemma
\begin{lemma}
Let $Y$ be a compact manifold without boundary. Suppose that $Q\in \Psi^0(Y,\complex ^k\otimes \complex ^k)$ with $Q^2-Q=0$, and $B\in \Psi^m(Y;\complex^k\otimes \complex^n)$ with symbols $q$ and $b$ respectively. Then,
\begin{enumerate}
\item if $b(y,\eta)$ restricted to $q(y,\eta)\complex^k$ is surjective for all $(y,\eta)\in T^*(Y)$, then one can find $S\in \Psi^{-m}(Y;\complex^n\otimes \complex^k)$ such that 
$$BS=I_n+O_{\Psi}(h),\quad QS=S$$
\item if $b(y,\eta)$ restricted to $q(y,\eta)\complex^k$ is injective for all $(y,\eta)\in T^*(Y)$, then one can find $S'\in \Psi^{-m}(Y;\complex^n\otimes \complex^k)$ and $S''\in \Psi^0(Y;\complex^k\otimes \complex^k)$ such that 
$$S'B+S''=I_k+O_{\Psi}(h),\quad S''Q=0$$
\item if $b(y,\eta)$ restricted to $q(y,\eta)\complex^k$ is bijective, then $S$, $S'$, $S''$ are uniquely determined mod $O_{\Psi}(h)$ and $S'=S+O_{\Psi}(h)$. 
\end{enumerate}
\end{lemma}
\noindent{\bf Remark:} $I_n$ denotes $i(x,hD)$ where $i$ is the $n\times n$ identity matrix.

\noindent \begin{proof}
To prove the first claim, observe that $bq$ is surjective. Hence there exists a right inverse $c$. Thus, $bqc=i$ and, letting $C=c(x,hD)$, $BQC=I_n+O_{\Psi}(h)$ and $Q(QC)=QC$. Hence, the first claim follows from letting $S=QC$. 

To prove the second claim, observe that $b\oplus (i-q)$ is injective and hence has a left inverse $(t',t'')$. Thus, letting $S'=t'(x,hD)$ and $S''=t''(x,hD)$, we have the claim.

To prove the third claim, just observe that 
$$S'=S'BS+O_{\Psi}(h)=S'BS+S''(S-QS)+O_{\Psi}(h)=S'BS+S''S+O_{\Psi}(h)=S+O_{\Psi}(h).$$
\end{proof}

 $Q$ and $B$ satisfy the hypotheses of the part (iii) of the previous lemma by \cite[Theorem 19.5.3]{HOV3}. Therefore, using this on $Q$ and $B$ from above, we obtain $S$ and $S''$ as described.

As in \cite[Chapter 20]{HOV3}, we use the properties of $S$ and $S''$ to see that 
$$R:(f,g)\mapsto (I+TP^cS''\gamma)Tf^0+TP^cSg$$
has 
$$PR(f,g)= f+ O(h)(f+g)\quad B_j\gamma R(f,g)=g_j+O(h)(f+g)$$
and 
$$R(Pu,\gamma u)=u+O(h)u$$
provided that $R$ is bounded 
$$R:H_h^{s-m}(\Omega) \oplus_{j=0}^{m-1}H_h^{s-m_j-\recip{2}}(\partial\Omega)\to H_h^{s}(\Omega).$$
Hence $R$ is a candidate for an approximate inverse modulo $O(h)$ errors. 

Therefore, in order to show that $R$ is both an approximate left and right inverse for (\ref{eqn:BVP}), all we need is the following lemma which follows from a rescaling of \cite[Proposition 20.1.6]{HOV3} along with the fact that $PT=I$ with no remainder. We include the proof here for convenience.
\begin{lemma}
If $s\geq m$ and $f\in C^\infty(\overline{\Omega})$, then 
\begin{equation}
\label{eqn:part1}\|T(f^0)\|_{H_h^s}\leq C\|f\|_{H_h^{s-m}}.\end{equation}
If $U=(U_0,...,U_{m-1})\in C^\infty(\partial \Omega)$, we have for any $s$
\begin{equation}\label{eqn:part2}\|TP^cU\|_{H_h^s}\leq C\sum_0^{m-1}\|U_j\|_{H_h^{s-j-\recip{2}}}.\end{equation}
\end{lemma}
\begin{proof}
It suffices to prove (\ref{eqn:part1}) when $f$ has support in a compact subset $K$ of a local coordinate patch $Y\times [0,1)$ at the boundary, $Y\subset \re^{d-1}$. Let $\chi\in C_0^\infty(K\times [0,1))$ have $\chi\equiv 1$ in a neighborhood of $K$, let $k\geq s-m$. Then, with the notation $\xi=(\xi_1,\xi ')$, and
$$\|u\|^2_{H_h^{(m,s)}}:=(2\pi h)^{-d}\int |\mcal{F}_h(u)(\xi)|^2\la \xi\ra^{2m}\la \xi '\ra^{2s}d\xi,$$
where $\mcal{F}_h$ is the semiclassical Fourier transform,
we have 
$$\|f^0\|_{H_h^{(-m+s-k,k)}}\leq \|f^0\|_{H_h^{(0,s-m)}}\leq \|f\|_{H_h^{s-m}}.$$
We can write 
$$(hD)^\alpha \chi T=\sum_{|\beta|\leq |\alpha|}T_\beta (hD)^\beta$$
where $T_\beta$ is a pseudodifferential operator of order $-m$ and $\beta_1=0$ if $\alpha_1=0$. Then, since $\|(hD)^\beta f^0\|_{H_h^{s-m-k}}\leq \|f^0\|_{H_h^{(s-m-k,k)}}$ if $|\beta|\leq k$ and $\beta_1=0$, we have
$$\|(hD)^\alpha\chi Tf^0\|_{H_h^{s-k}}\leq C\|f^0\|_{H_h^{(s-m-k,k)}}$$
if $\alpha_1=0$ and $|\alpha|\leq k$. Thus, 
$$\|\chi Tf^0\|_{H_h^{(s-k,k)}}\leq C\|f\|_{H_h^{(s-m)}}.$$
But, $PTf^0=f$ in $\Omega^o$. Hence, by \cite[Theorem B.2.9]{HOV3}, or rather a rescaling of its proof, we have for $\psi\in C_0^\infty$ with $\psi\equiv 1 $ in a neighborhood of $K$ and $\chi\equiv 1$ in a neighborhood of supp $\psi$, 
\begin{equation}
\label{eqn:commutex}\|\psi Tf^0\|_{H_h^s}\leq C(\|f\|_{H_h^{s-m}}+\|\chi Tf\|_{H_h^{(s-k,k)}})\leq C\|f\|_{H_h^{s-m}}.
\end{equation}
Since $\|(1-\psi)Tf^0\|_{H_h^s}=O(h^\infty)\|f^0\|_{H_h^{s-m}}$, we have proved (\ref{eqn:part1}). 

Now, to prove (\ref{eqn:part2}), we may assume that supp $U\subset K$. We have, 
$$P^cU=\sum_{0}^{m-1}v_j\otimes (hD_1)^j\delta,\quad v_j=\sum_{j+l<m}P_{j+l+1}U_lhi^{-1}.$$
Then, since $P_{j+l+1}$ is order $m-j-l-1$, we have 
$$\sum \|v_j\|_{H_h^{s-m+j+\recip{2}}}\leq Ch\sum\|U_j\|_{H_h^{s-j-\recip{2}}}.$$
The semiclassical Fourier transform of $v_j\otimes (hD_1)^j\delta$ is $\mcal{F}_h(v_j)(\xi ')\xi_1^j$, and, when $j<m$,
$$h^{-1}\int \xi_1^{2j}(1+|\xi '|^2+\xi_1^2)^{-m}d\xi_1\leq Ch^{-1}(1+|\xi '|^2)^{j-m+\recip{2}}.$$
Thus, if $j<m$ and $k$ is an integer, $k\geq \max(s,0).$
$$\|v_j\otimes (hD_1)^j\delta\|_{H_h^{(-m+s-k,k)}}\leq \|v_j\otimes (hD_1)^j\|_{H_h^{(-m,s)}}\leq Ch^{-1}\|v_j\|_{H_h^{s-m+j+\recip{2}}}.$$
Putting these together, we have 
$$\|TP^cU\|_{H_h^{(s-k,k)}}\leq C\sum_{0}^{m-1}\|U_j\|_{H_h^{s-j-\recip{2}}}.$$
Then, because $T$ is continuous from $H_h^{(t-m,k)}$ to $H_h^{(t,k)}$ for $k\geq 0$, we can commute $T$ with $x'$ derivatives. Then, by observing that $PTP^cU=P^cU$ and proceeding as in \eqref{eqn:commutex} we can improve this estimate to \eqref{eqn:part2}. 
\end{proof}

Proposition \ref{prop:Calderon} now follows from the fact that $R$ is an inverse for \eqref{eqn:BVP} modulo $O(h)$ errors and an application of a standard Neumann series argument that can be found, for example, in \cite[Theorem C.3]{EZB}.

\section{Construction of Quasimodes Via Boundary WKB method}
\label{sec:QuasiConstruct}
In this section, we will prove part 3 of Theorem \ref{thm:PseudoSpec}. Moreover, we do not assume that $X$ is constant in the construction. In particular, we show
\begin{prop}
\label{prop:constructQuasi}
Let $X\in C^\infty(\re^d;\re^d)$ be a vector field and $\partial \Omega_+$ defined as in \eqref{eqn:OmegaBound}. Then, for each $x_0\in \partial \Omega_+$, let $\nu(x_0)$ be the outward unit normal. Then, if $d\geq 2$, for each
$$z\in \{\zeta\in \complex:|X(x_0)|^{-2}(\Im \zeta)^2<\Re \zeta\}\setminus \{\la X(x_0),\nu(x_0)\ra^2/4\},$$ 
there exists $u\in C^\infty(\overline{\Omega})$ such that $u$ is a quasimode for $z$ with $\ES_h(u)=x_0$. Moreover, if $\partial \Omega$ and $X$ are real analytic near $x_0$, then $P_zu=O_{L^2}(e^{-\delta/h}).$	

If $d=1$, then for each $x_0\in \partial\Omega_+$ and each 
$$z\in \{\zeta\in \complex:|X(x_0)|^{-2}(\Im \zeta)^2<\Re \zeta\}\setminus\{\zeta\in \complex:\Im \zeta=0\text{ and }\Re \zeta\geq \la X(x_0),\nu(x_0)\ra^2/4\},$$
there exists $u\in C^\infty(\overline{\Omega})$ with the same properties as above.
\end{prop} 

\vbox{
\noindent{\bf Remarks: }
\begin{enumerate}
\item We demonstrate the construction in dimension $d\geq 2$. The additional restriction in $d=1$ comes from the fact that $\partial \Omega$ is a discrete set of points and hence functions on $\partial\Omega$ are determined by their values at these points. In particular, since we cannot choose $d\phi_0$ for $\phi_0$ in equation \eqref{eqn:phase}, we must restrict the $z$ for which we make the construction.
\item In fact we will also show that $P_zu=O_{L^\infty}(h^\infty)$ in the smooth case and $O_{L^\infty}(e^{-\delta/h})$ in the analytic case.
\end{enumerate}
}

We wish to construct a solution to (\ref{eqn:generic}) that concentrates at a point in $\partial \Omega_+$. Let $x_0\in \partial \Omega_+$ and assume for simplicity that $|X(x_0)|=1$ and without loss that $X(x_0)=e_1$. We also assume $\Im z\neq \frac{\nu_1^2}{4}$ for technical reasons. To accomplish the construction, we postulate that $u$ has the form 
\begin{equation}
\label{eqn:ampExpand}
u=\chi v,\quad v=(a e^{i\varphi_{1,N}/h}-be^{i\varphi_{2,N}/h}),\quad 
a=\sum_{n=0}^N a_n h^n,\quad b=\sum_{n=0}^N b_nh^n.
\end{equation}
Then, let $\Gamma\subset \partial \Omega$ be a small neighborhood of $x_0$ to be determined later and $U\subset \overline{\Omega}$ be a small neighborhood of $\Gamma$. We solve for $a$, $b$, $\varphi_{1,N}$, and $\varphi_{2,N}$ such that
\begin{equation}\begin{cases}
\label{eqn:quasimode} P_zv=O_{L^2}(h^{N+2})&\text{ in }U,\\
v|_{\Gamma}=0.\end{cases}\end{equation}

More precisely, we find two distinct solutions, $\varphi_{1,N}$ and $\varphi_{2,N}$ to 
\begin{equation}
\label{eqn:phase}
\begin{cases} p_z(x,d \varphi_{i,N})=O(|x-x_0|^{2N+4})&\text{ in }U,\\
\varphi_{i,N}|_{\Gamma }=\phi_0,\end{cases}
\end{equation}
with $\phi_0(x_0)=0$,

\begin{equation}
\label{eqn:phiChoice}
d\phi_0(x_0)=\lambda\left\{\begin{array}{cl} e_1'&|\nu_1|<1\text{ or }d=1\\
 e_2&\nu_1=1,\text{ }d\geq 2\end{array}\right\} \in T^*\Gamma,
\end{equation}

\noindent $\Im d\phi_0=0$ and $\Im d^2\phi_0(x_0)>0$, where $e_1'=(e_1-\la e_1,\nu(x_0)\ra\nu )/X'$, $X'=\sqrt{1-\nu_1^2}$, $\nu_1=\la e_1,\nu(x_0)\ra,$ and $\lambda$ will be chosen later (see \eqref{eqn:lambda}). 

\noindent {\bf Remark:} We choose $\phi_0$ in this way to get localization along the boundary.

In addition, we solve the transport equations
\begin{equation}
\label{eqn:amplitude}
\begin{cases}-i\Delta \varphi_{i,N} \psi_n+2i\la \partial \varphi_{i,N}, \partial \psi_n\ra +\la X,\partial \psi_n\ra=\Delta \psi_{n-1}+O(|x-x_0|^{2N+4})&\text{ in }U,\\
\psi_0|_{\Gamma}=1\quad \psi_n|_{\Gamma}=0\text{ for }n>0,\end{cases}.
\end{equation}
\eqref{eqn:phase} has two solutions ($\varphi_{1,N}$ and $\varphi_{2,N}$) and we set $a_n:=\psi_n$ when using $\varphi_{1,N}$ in \eqref{eqn:amplitude} and $b_n:=\psi_n$ when using $\varphi_{2,N}$. 

First, we consider (\ref{eqn:phase}). To solve this equation, we construct a complex Lagrangian submanifold as in \cite[Theorem 1.2']{zworski04}. Note that with the choice of $d\phi_0(x_0)$ as in \eqref{eqn:phiChoice}, $$d \varphi(x_0)=(\alpha+i\beta)\nu(x_0)+d\phi_0(x_0)$$
for some $\alpha,\beta\in \re$. This gives rise to
$$\alpha^2-\beta^2-\beta\nu_1-\Re z+\lambda^2=0$$
and
$$\alpha(2\beta+\nu_1)-\Im z+X'\lambda=0.$$

Hence
$$(\beta^2+\beta\nu_1+\Re z-\lambda^2)(2\beta+\nu_1)^2=(\Im z- X'\lambda)^2.$$
Letting $c=\beta+\frac{\nu_1}{2}$, we have  
$$c^4+\left(\Re z-\lambda^2-\frac{\nu_1^2}{4}\right)c^2-\recip{4}(\Im z-X'\lambda)^2=0.$$
Which gives 
$$c=\pm \recip{\sqrt{8}}\left(\sqrt{\nu_1^2-4(\Re z-\lambda^2)+ \sqrt{(\nu_1^2-4(\Re z-\lambda^2))^2+16(\Im z-X'\lambda)^2}}\right)$$
(note that we take the positive root inside so that the result is real).

In order to complete the construction, we need $g<0$. That is, we require
\begin{eqnarray}
\frac{\nu_1}{2}>|c|&\Leftrightarrow& 2\nu_1^2>\nu_1^2-4(\Re z-\lambda^2)+ \sqrt{(\nu_1^2-4(\Re z-\lambda^2))^2+16(\Im z-X'\lambda)^2}\nn
&\Leftrightarrow&(\nu_1^2+4(\Re z-\lambda^2))^2>(\nu_1^2-4(\Re z-\lambda^2))^2+16(\Im z-X'\lambda)^2\nn 
&\Leftrightarrow&(\Re z-\lambda^2)\text{ }\nu_1^2>(\Im z-X'\lambda)^2\nn 
&\Leftrightarrow&\Re z>\frac{(\Im z )^2(1-X'^2)+(\lambda-\Im z X')^2}{\nu_1^2}=(\Im z)^2+\frac{(\lambda-\Im zX')^2}{\nu_1^2}.\nonumber
\end{eqnarray}

We also need $|c|>0$ so that $\varphi_{1,N}$ and $\varphi_{2,N}$ are distinct. Thus, letting $|a|<1$, we choose 
\begin{equation}
\label{eqn:lambda}\lambda =\begin{cases}\Im zX'&\Im z\neq 0,\\
\nu_1\sqrt{\Re z}a&\Im z=0, \nu_1\neq 1,\\
0&\Im z=0, \Re z< \recip{4}, \nu_1=1,\\
\sqrt{\Re z-\recip{8}}&\Im z=0, \Re z\geq \recip{4}, \nu_1=1.\end{cases}\end{equation}
\noindent{\bf Remark:} In dimension 1, we are forced to choose $\lambda=0$, however, in dimension 1, $\nu_1=1$ and $X'=0$, so we have $|c|>0$ when $\Re z<\recip{4}=\frac{\nu_1^2}{4}$ or $\Im z\neq 0$.

With this choice for $\lambda$, we have $0<|c|<\frac{\nu_1}{2}$ if and only if $(\Im z)^2<\Re z$. Hence, $v$ decays exponentially in the $-x_1$ direction if $(\Im z)^2<\Re z$. This decay allows us to localize our construction near the boundary.

\noindent {\bf Remark:} $(\Im z)^2<\Re z$, corresponds precisely with (\ref{eqn:ellipticRestriction}). 

Now that we have $f$, $g$, and $\lambda$, we need to solve (\ref{eqn:phase}) on the rest of $\Gamma$ and in the interior of $U$.

\subsection{Analytic Case}
We first assume that $\Gamma$ and $p$ are real analytic and solve the equations exactly. Let $\gamma$ be the coordinate change to semigeodesic coordinates for $\Gamma.$ Extend $\gamma$ analytically to a neighborhood of $\Gamma$ in $\complex^{d-1}$. Then, define $\kappa$, the lift of $\gamma$, by $(z,\eta)\mapsto (\gamma(z),(\partial \gamma^{-1})^T\eta).$ Next, choose $\phi_1(y)$ real analytic in a neighborhood of $0\in \re^{d-1}$ with $\phi_1(0)=0$, $d\phi_1(0)$ as in \eqref{eqn:phiChoice}, and $\Im d^2\phi_1(0)>0$. Then, extend $\phi_1$ to $y$ in a neighborhood of the origin in $\complex^{d-1}$. Next, let 
$$\Lambda_0:=\{(0,y,\xi_1(y),d_y\phi_1(y)):\kappa_*p_z(0,y,\xi_1,d_y\phi_1(y))=0,(0,y)\in \gamma(\Gamma)\}$$ 
where $\xi_1(y)$ is well defined since $\xi_1(0)=f+ig$ and, for $z\neq \frac{\nu_1^2}{4}$, $\partial_{\xi_1}\kappa_*p_z(x_0,f+ig,d\phi(x_0) )\neq 0$.
Observe also that $\Lambda_0$ is isotropic with respect to the complex symplectic form.

Finally, let $\Phi_t$ be the complex flow of $\kappa_*p_z$ which exists by the Cauchy-Kovalevskaya Theorem (\cite[Section 4.6]{E}). Then, $$\Lambda:=\cup_{|t|<\e}\Phi_t(\Lambda_0)$$
is Lagrangian. Hence it has a generating function $\tilde{\varphi}$ such that $\varphi=\tilde{\varphi}\composed \gamma$ solves \eqref{eqn:phase} and has $\varphi |_\Lambda=\phi_0:=\phi_1\composed\gamma$. Therefore, there exist $\varphi_{1,N}\neq \varphi_{2,N}$, solutions to (\ref{eqn:phase}).

\noindent {\bf{Remark:}} Note that the two distinct solutions $\varphi_{1,N}$ and $\varphi_{2,N}$ come from the two solutions to $\xi_1(x_0)$. 

Next, we solve (\ref{eqn:amplitude}). To do this, note that $\varphi_{1,N}$ and $\varphi_{2,N}$ from above are analytic. Hence, since $\Gamma$ and (\ref{eqn:amplitude}) are analytic, we may apply the Cauchy-Kovalevskaya Theorem as above to find $a_n$ and $b_n$.

If $\Gamma$ and $p$ are analytic, it is classical \cite[Theorem 9.3]{SjoMicro} that the solutions $a_n$ and $b_n$ have $\max(|a_n|,|b_n|)\leq C^nn^n.$ This will be used below to show that the error contributed by truncation at $N=1/Ch$ is exponential.

\subsection{Smooth Case} 
Suppose that $\Gamma$ and $p$ are not analytic. Then, let $\gamma$ be the coordinate change to semigeodesic coordinates for $\Gamma$. Define the lift $\kappa$ of $\gamma$ and choose $\phi_1$ as above. We now solve the equations (\ref{eqn:phase}) with $O(|x-x_0|^{2N+4})$ error. First, write 
$$\kappa_*p_z(x,\xi)=p_1(x,\xi)+O(|x-x_0|^{2N+4})p'(\xi),$$ 
where $p_1$ is the Taylor polynomial for $\kappa_*p_z$ to order $2N+4$. Next, apply the construction for analytic $p_z$ from above to solve 
$$\begin{cases}
p_1(x,d \theta)=0,\\
\theta|_{\gamma(\Gamma)}=\phi_1.
\end{cases}$$
Then, observe that 
$$\kappa_*p_z(x,d\theta)=O(|x-x_0|^{2N+4})p'(d\theta)=O(|x-x_0|^{2N+4}).$$
Hence, we have that $\varphi_{i,N}=\theta\composed \gamma$ solves \eqref{eqn:phase}. 

Now, using the solution $\varphi_{i,N}$ just obtained, we solve the amplitude equations $(\ref{eqn:amplitude})$ with $O(|x-x_0|^{2N+4})$ errors. As with the phase, we start by changing to semigeodesic coordinates. Write the equation for $\psi_n$ in the new coordinates as 
$$\begin{cases}
\la \rho, \partial \psi'_n\ra- \zeta \psi'_n= f,\\
\psi'_0|_{\gamma(\Gamma)}=1\quad \psi'_n|_{\gamma(\Gamma)}=0\text{ for }n>0.
\end{cases}$$
Then, writing the Taylor polynomials to order $2N+4$ for $\rho,\zeta,$ and $f$ as $\rho_1,\zeta_1$, and $f_1$ respectively, we solve 
$$\begin{cases}
\la \rho_1, \partial \psi'_n\ra- \zeta_1\psi'_n= f_1,\\
\psi'_0|_{\gamma(\Gamma)}=1\quad \psi'_n|_{\gamma(\Gamma)}=0\text{ for }n>0,
\end{cases}$$
using the analytic construction above. Then, just as in the solution of \eqref{eqn:phase}, $\psi_n:=\psi'_n\composed\gamma$ solves \eqref{eqn:amplitude}.

\noindent{\bf Remark:} We are actually solving for the formal power series of $\varphi_{i,N}$, $a_n$ and $b_n$. 

\subsection{Completion of the construction}
Let $V\Subset U$ be a neighborhood of $\Gamma$. Then, let $\chi\in C^\infty(\overline{\Omega})$ with $\chi\equiv 1$ on $V$ and $\chi\equiv 0$ on $\overline{\Omega}\setminus U$. For convenience, we make another change of coordinates so $x_0\mapsto 0$ and that $\nu(x_0)=e_1$. Then, $\Im \partial_{x_1}\varphi_{i,N}(0)<0$ for $i=1,2$, and $\Im d^2\phi(0)>0$. Together, these imply that $\Im \varphi_{i,N} >0$ on supp $\partial \chi\cap \Gamma$. Hence, we have for $U$ small enough but independent of $h$,
$$|P_zu|=\left|\chi P_zv+[P_z,\chi]v\right|\leq O(|x|^{2N+4})\left(\left|e^{i\varphi_{1,   N}/h}\right|+\left|e^{i\varphi_{2,N}/h}\right|\right)+O(h^{N+2})+O(e^{-\e/h}).$$  
Now, observe that $a|_\Gamma =b|_\Gamma$, $\varphi_{1,N}|_\Gamma =\varphi_{2,N}|_\Gamma=\phi_0$ and 
$$\varphi_{i.N}=\phi_0(x')+c_ix_1+O(x_1|x'|)+O(x_1^2),$$
 $c_1\neq c_2$. Hence, since $\Im \phi_0(x')\geq c |x'|^2$, and $\Im c_i<0$ 
$$O(|x|^{2N+4})(e^{i\varphi_{1,N}/h}+e^{i\varphi_{2,N}/h})=O_{L^2\cap L^\infty }(h^{N+2})$$
and $v$ solves \eqref{eqn:quasimode}.

Note also that if $\Gamma$ and $p$ are analytic, then the equations \eqref{eqn:phase} and \eqref{eqn:amplitude} can be solved exactly with $\max(|a_n|,|b_n|)<C^nn^n$. Hence, truncating the sums \eqref{eqn:ampExpand} at $N=1/eCh$, we have 
\begin{align*}
|P_zu|=\left|\chi P_zv+[P_z,\chi]v\right|\leq C^NN^Nh^{N+1} + O(e^{-\e/h})&=(CN)^N(CN)^{-N-1}e^{-N-1}+O(e^{-\e/h})\\&=e^{-c/h}+O(e^{-\e/h})=O(e^{-\e/h}).\end{align*}

Our last task is to show that $\|u\|_{L^2}\geq Ch^{\frac{d+3}{4}}$. To see this, we calculate, shrinking $U$ and $V$ if necessary, and letting $u=\chi v$,
\begin{eqnarray}
\|u\|_{L^2}^2&=&\int_V\left|(1+O(x))\left(e^{i\varphi_{1,N}/h}-e^{i\varphi_{2,N}/h}\right)\right|^2\nn 
&\geq &c\int_{\gamma (V)}\left|e^{i(\phi_0(x')+O(|x'|x_1)+O(x_1^2))/h}\left(e^{ic_1x_1/h}-e^{ic_2x_1/h}\right)\right|^2dx_1dx'\nn 
&\geq &c\int_{\gamma (V)}e^{-c|x'|^2/h}e^{-c_0x_1/h}x_1\left(1+O(e^{-\delta x_1/h})\right)dx_1dx'\geq ch^{\frac{d+3}{2}}\nonumber
\end{eqnarray}

\noindent{\bf Remark:} By the same argument $\|u\|_{L^2}^2\leq Ch^{\frac{d+3}{2}}$.

To finish the construction of $u$, we simply rescale $u$ so that it has $\|u\|_{L^2}=1$ and invoke Borel's Theorem (see, for example \cite[Theorem 4.15]{EZB}). 

\section{Propagation of Semiclassical Wavefront Sets}
\label{sec:propagation}

We first examine the case where $P=hD_{x_1}+ihD_{x_2}=hD_{\bar{z}}$ (here, we identify $\re^2$ with $\complex$). We make the following definition in the spirit of Duistermaat and H\"{o}rmander \cite[Section 7]{FIO2}
\begin{defin} 
$$s^0_u(x):=\sup\{t\in \re:\text{ there exists } U \text{ a neighborhood of }x\text{ such that } \|h^{-t}u\|_{L^2(U)}=O(1)\}$$
\end{defin}
We will need the following lemma 

\begin{lemma}
\label{lem:superHarmx}
Let $u\in H_h^{1}$ with $hD_{\bar{z}}u=f$. Then, $s^*=\min(s^0_u,s)$ is superharmonic if $s$ is superharmonic and $s^{0}_f-1\geq s$. 
\end{lemma}
\begin{proof}
Let $q(z)$ be harmonic function in $\complex$ such that $s^*(z,0)\geq q(z)$ for $|z|=r$. (Here we have written $x\in \re^d$ as $x=(z,x_3,...,x_d)$, identifying $\complex$ with $\re^2$.) Then, we need to show that $s^*(z)\geq q(z)$ for $|z|\leq r$. The fact that $s(z,0)\geq q(z)$ for $|z|\leq r$ follows from the superharmonicity of $s$. Therefore, we only need to show the inequality for $s^0_u$.

To do this, let $\chi_1\in C_0^\infty(\re^{d-2})$ have support in a small neighborhood of $0$, $\chi_1(0)=1$, and $\chi_2\in C_0^\infty (\complex)$ be 1 for $|z|\leq r$ and $0$ outside a neighborhood so small that
$$\|h^{-q(z)}u\|_{L^2(\text{supp }\chi_1\partial \chi_2)}=O(1).$$
This is possible since $s^*\geq q(z)$ for $|z|=r$ implies $s^0_u\geq q(z)$ for $|z|=r$. Then, define $v:=\chi_1\chi_2u$. We have
$$hD_{\bar{z}}v=\chi_1\chi_2f+\chi_1uhD_{\bar{z}}(\chi_2)=:g$$
with $s^{0}_g(z)\geq q(z)+1$. 

Next, let $F(z)$ be analytic with $\Re F(z)=q(z)$ and define $Q(z):=h^{-F(z)}$. Using this, we have 
$D_{\bar{z}}(Qv)=h^{-1}Qg$ and, since $s^{0}_g\geq q+1$, $h^{-1}Qg=O_{L^2}(1).$ Then, applying $\partial_z$, we have 
$$\Delta_{x_1,x_2} (Qv)=CD_zD_{\bar{z}}(Qv)=Ch^{-1}D_z(Qg)$$
and hence, shrinking $r$ if necessary (note that this is valid since superharmonicity is a local property), $Qv$ solves,
$$\begin{cases}
\Delta_{x_1,x_2}(Qv)=Ch^{-1}D_z(Qg)&\text{in }B(0,1)\subset \re^2,\\
Qv=0&\text{in }\partial B(0,1).
\end{cases}$$
Therefore, by the estimate for $u$ with $u|_{B(0,1)}=0$,
$$\|u\|_{L^2}\leq C\|\Delta u\|_{H^{-1}},$$
 we have
$$\|Qv\|_{L^2_{x_1,x_2}}\leq C\|h^{-1}D_zQg\|_{H^{-1}_{x_1,x_2}}.$$
But, since $h^{-1}Qg=O_{L^2}(1)$, $\|h^{-1}D_zQg\|_{H^{-1}_{x_1,x_2}}=O(1)$ for almost every $x'\in $ supp $\chi_1$, the same is true for $Qv$. 
Thus, since $v\equiv u$ in a neighborhood of $0$, and $|Q|=|h^{-q(z)}|$, $s^0_u\geq q(z)$ for $|z|\leq r$. 
\end{proof}

\begin{defin}
\label{defn:quantizeTrans}
We say that an operator $T$ quantizes $\kappa$ if $T:L^2\to L^2$ and for all $a\in S^m$, we have 
$$T^{-1}a^w(x,hD)T=b^w(x,hD)$$
for a symbol $b\in S^m$ satisfying 
$$b|_{U_0}:=\kappa^*(a|_{U_1})+O_{S^{m-1}}(h).$$
\end{defin}

To convert from $P$ as in \eqref{eqn:P} to the case of $P=hD_{\bar{z}}$ we need the following lemma similar to \cite[Theorem 12.6]{EZB} which we include for completeness. 
\begin{lemma}
\label{lem:canonTrans}
Suppose $P=p^w$ and $p$ has 
$$\sum_{k=0}^\infty h^kp_k$$
with $p_0(0,0)=0$ and $\{\Re p, \Im p\}=0$ with $\partial \Re p$ and $\partial \Im p$ linearly independent. Then there exists a local canonical transformation $\kappa$ defined near $(0,0)$ such that 
$$\kappa^*p_0=\xi_1+i\xi_2$$
and an operator $T:L^2\to L^2$ quantizing $\kappa$ in the sense of Definition \ref{defn:quantizeTrans} such that $T^{-1}$ exists microlocally near $((0,0),(0,0))$ and $$TPT^{-1}=hD_{x_1}+ihD_{x_2}\quad \text{microlocally near }((0,0),(0,0)).$$
\end{lemma}
\begin{proof}
Let $q_1=\Re p_0$ and $q_2=\Im p_0$. Then, by a variant of Darboux's Theorem (see, e.g. \cite[Theorem 12.1]{EZB}), there exists $\kappa$ a symplectomorphism, locally defined near $(0,0)$, such that $\kappa(0,0)=(0,0)$ and 
$$\kappa^*q_1=\xi_1\quad \kappa^*q_2=\xi_2.$$
Then, by \cite[Theorem 11.6]{EZB} shrinking the domain of definition for $\kappa$ if necessary, there exists a unitary $T_0$ quantizing $\kappa$ such that 
$$T_0PT_0^{-1}=hD_{x_1}+ihD_{x_2}+E\quad\text{ microlocally near }(0,0),$$
where $E=e^w$ for $e\in hS$ .

Next, we find $a\in S$ elliptic at $(0,0)$ such that 
$$hD_{x_1}+ihD_{x_2}+E=A(hD_{x_1}+ihD_{x_2})A^{-1}\text{ microlocally near }(0,0),$$
where $A=a^w$
i.e.
$$[hD_{x_1}+ihD_{x_2},A]+EA=0\text{ microlocally near }(0,0).$$

Since $P=p_0^w+hp_1^w+...$, we have $E=e^w$ for $e=he_1+h^2e_2+...$. We use the Cauchy formula to solve the equation 
$$\recip{i}\{\xi_1+i\xi_2,a_0\}+e_1a_0=0$$
near $(0,0)$ for $a_0\in S$ with $a_0(0,0)\neq 0$. Then, defining $A_0:=a_0^w$, we have 
$$[hD_{x_1}+ihD_{x_2},A_0]+EA_0=r_0^w$$
for $r_0\in h^2S$.
To complete the proof, we proceed inductively to obtain $A_k=a_k^w$ for $a_k\in S$, solving
$$[hD_{x_1}+ihD_{x_2},A_0+...+h^NA_N]=E(A_0+...+h^NA_N)=r_N^w$$
where $r_N\in h^{N+2}S$, using the Cauchy formula at each stage. Then, we invoke Borel's Theorem (see, for example \cite[Theorem 4.15]{EZB}) to find $A$ and let $T=A^{-1}T_0$. 
\end{proof}

Now, define 
\begin{defin}
\label{def:WFhmeasureTemp}
\begin{multline*}
S^{0}_u(x,\xi):=\sup \{t\in \re :\exists \text{ }U,V\subset \re^d\text{ open } x\in U,\text{ }\xi\in V\text{ s.t. }\\
 \forall \text{ }\chi_1\in C_0^\infty(U),\chi_2\in C_0^\infty(V), h^{-t}\chi_1(x)\chi_2(hD)u=O_{L^2}(1) \}
\end{multline*}
\end{defin}

\noindent In Lemma \ref{lem:relationtoWFH} we prove that Definition \ref{def:WFhmeasureTemp} is equivalent to the following
\begin{defin}
\label{def:WFhmeasure}
\begin{multline*}
S_u(x,\xi):=\sup \{t\in \re :\text{ there exists }U\subset\re^{2d}\text{ open}\text{, } (x,\xi)\in U,\\
\text{ s.t. } \forall \text{ }\chi\in C_0^\infty(U)\text{ } h^{-t}\chi^w(x,hD)u=O_{L^2}(1) \}.$$
\end{multline*}
\end{defin}

\noindent The proof follows \cite[Theorem 8.13]{EZB}, but we reproduce it in this setting for the convenience of the reader.
\begin{lemma}
\label{lem:relationtoWFH}
Suppose that there exist $U$ and $V$ as in Definition \ref{def:WFhmeasureTemp}. Then, there exists $W$ open, $(x,\xi)\in W$ such that for $\chi\in C_0^\infty(W)$, $\chi^wu=O_{L^2}(h^{t}).$
\end{lemma}
\begin{proof}
Let $a=\chi_1(x)\chi_2(\xi)$ as in Definition \ref{def:WFhmeasureTemp}. Then, there exists $\chi\in C_0^\infty(\re^{2d})$ supported near $(x_0,\xi_0)$ such that 
$$|\chi(x,\xi)(a(x,\xi)-a(x_0,\xi_0))+a(x_0,\xi_0)|\geq \gamma >0.$$
Hence, by \cite[Theorem 4.29]{EZB}, there exists $c\in S$ such that for $h$ small enough, 
$$c^w(\chi^wa^w+a(x_0,\xi_0)(I-\chi^w))=I.$$

Next, observe that 
$$b^wu=b^wc^w\chi^w a^wu+a(x_0,\xi_0)b^wc^w(I-\chi^w)u.$$
Now, the first term on the right is $O_{L^2}(h^{t})$ since $a^wu=O_{L^2}(h^t)$. Also, if the support of $b$ is sufficiently near $(x_0,\xi_0)$, supp $b\cap$ supp $(1-\chi)=\emptyset$ and hence the second term is $O_{L^2}(h^\infty)$. This proves the claim.
\end{proof}

\noindent{\bf Remark:} Note that $S_u(x,\xi)=\infty$ if and only if $(x,\xi)\notin \WF_h(u).$

Lemma \ref{lem:relationtoWFH} shows that $S_u(x,\xi)=S_u^0(x,\xi)$. It will be convenient to use both of these definitions in the proof of the following proposition.

\begin{prop}
\label{prop:prop}
Let $u\in H_h^m$, $p(x,hD)u=f$ with $p\in S^m$ and let $S_f\geq s+1$, $O \subset N$ where $$N:=\left\{(x,\xi) \in  T^*\re^d: p(x,\xi)=0,\text{ }\{p, \bar{p}\}=0,\text{ } H_{\Re p}\text{ and }H_{\Im p} \text{ are independent}\right\}.$$ Then, it follows that $\min (S_u, s)$ is superharmonic in $O$ if $s$ is superharmonic in $O$, and that $\min(S_u-s,0)$ is superharmonic in $O$ if $s$ is subharmonic in $O$ (with respect to $H_p$). In particular, $S_u$ is superharmonic in $O$ if $O\cap \WF_h(f)=\emptyset$. 
\end{prop}
\begin{proof}
First, we consider $hD_{\bar{z}}$. We prove that Lemma \ref{lem:superHarmx} remains valid with $s^0_u$ and $s^{0}_f$ replaced by $S^{0}_u$ and $S^{0}_f$.  Let $\chi\in C_0^\infty (\re^d)$. Then, $hD_{\bar{z}}(\chi(hD)u)=\chi(hD)f.$ Take $\chi_j$ with $\chi_j(\xi_0)=1$, vanishing outside $V_j$, $V_j\downarrow \{\xi_0\}$, and denote $u_j=\chi_j(hD)u$. Then we have $s^0_{u_j}(x)\uparrow S^{0}_u(x,\xi_0)$. So, the superharmonicity of 
$$\min (s^0_{u_j},s_j),\text{ where }s_j(x)=\inf_{\xi\in V_j}s(x,\xi),$$
gives that $\min (S^{0}_u,s)$ is superharmonic and proves the first part of the proposition for $hD_{\bar{z}}$. 

To prove the second, note that it is equivalent to the first if $s$ is harmonic. Thus, the second part follows if $s$ is the supremum of a family of harmonic functions. If $s\in C^2$ is strictly subharmonic, then $s(z,x',\xi)\geq q(z,x',\xi)$ in a neighborhood of $(w,x',\xi)$ with equality at $(w,x',\xi)$ when $q$ is the harmonic function 
$$q(z,x',\xi)=s(w,x',\xi)+\Re (2(z-w)\partial s(w,x',\xi)/\partial w+ (z-w)^2\partial^2s(w,x',\xi)\partial w^2).$$
Then, the local character of superharmonicity proves the second statement when $s$ is strictly subharmonic and the general case follows by approximation of $s$ with such functions. 

To pass from $hD_{\bar{z}}$ to $P$, we need the following (\cite[Lemma 7.2.3]{FIO2}).
\begin{lemma}
If $(x_0,\xi_0)\in N$ there exists $a\in S^{1-m}$ with $a(x_0,\xi_0)\neq 0$ such that $\{q,\bar{q}\}=0$ in a neighborhood of $(x_0,\xi_0)$ if $q=ap$. 
\end{lemma}
Now, by Lemma \ref{lem:canonTrans}, there exists $T$ microlocally quantizing $\kappa$ such that $\kappa^*(\Re ap)=\xi_1$ and $\kappa^*(\Im ap)=\xi_2$, so that
$$Ta^wPT^{-1}=hD_{x_1}+ihD_{x_2}$$
microlocally near $((0,0), (x_0,\xi_0))$. 
 Then, 
$$S_{Tu}\composed \kappa (x,\xi)=S_{u}(x,\xi),$$
for $(x,\xi)\in V$ a small neighborhood of $(x_0,\xi_0)$.
This follows from the fact that by Definition \ref{defn:quantizeTrans}, if $\chi\in C_0^\infty$, then
$$h^{-t}T^{-1} T\chi^wT^{-1} Tu=h^{-t}T^{-1}b^w Tu$$
where $b\in S$ and $b|_{U_0}=\kappa^*(\chi|_{U_1})+O_{S^{-1}}(h)$ and $T^{-1}$ is uniformly bounded on $L^2$.
But 
$$S^0_{Ta^wPu}\composed \kappa=S_{a^wf}=S_{f}=S^0_f\geq s+1.$$
Hence, Proposition \ref{prop:prop} follows from the case with $hD_{\bar{z}}$. 
\end{proof}

We need the following elementary lemma to prove Corollary \ref{cor:Propagation}. 
\begin{lemma}
\label{lem:WFSubsetCharSet}
Suppose $u$ solves (\ref{eqn:generic}) and $P_z$ has symbol $p_z(x,\xi)$. Then, 
$$\WF_h(u)\cap(\Omega^o\times\re^d)\subset p_z^{-1}(0)\cap(\Omega^o\times \re^d).$$
\end{lemma}
\begin{proof}
Let $x_0\in\Omega^o$ and $(x_0,\xi_0)\notin p_z^{-1}(0)$. Then, let $\chi_1\in C_0^\infty(\re^d)$ have support near $x_0$ and $\chi_2\in C_0^\infty(\re^d)$ have support near $\xi_0$. Then we have 
$$\chi_2^w\chi_1^wp_z^wu=O(h^\infty).$$
But, $\chi_2^w\chi_1^w=c^w$, $c\in S$ with $|c(x_0,\xi_0)|>0$. Similarly,
$c^wp_z^w=q^w$ for $q\in S$ with $|q(x_0,\xi_0)|>0$. Hence, $(x_0,\xi_0)\notin \WF_h(u)$.  
\end{proof}

Putting Proposition \ref{prop:prop} together with Lemma \ref{lem:WFSubsetCharSet}, we have the following corollary 
\begin{corol}
\label{cor:Propagation}
Let $P$ as in \eqref{eqn:P}, $\Re z>(\Im z)^2|X|^{-2}$, and $u\in H_h^2$ with $P_zu=O_{L^2}(h^\infty)$. Then, $\WF_h(u)\cap (\Omega^o\times \re^d)$ is invariant under the leaves generated by $H_{\Im p}=\la X,\partial_x\ra$ and $H_{\Re p}=2\la \xi, \partial_x\ra$. 
\end{corol}
\begin{proof}
By Lemma \ref{lem:WFSubsetCharSet}, 
$$\WF_h(u)\cap(\Omega^o\times \re^d)\subset p_z^{-1}(0)\cap(\Omega^o\times \re^d).$$ 
Also, $\{p,\bar{p}\}=0$, and, for $\Re z>(\Im z)^2|X|^{-2}$, $H_{\Re p}$ and $H_{\Im P}$ are independent on all of $p_z^{-1}(0).$
Now, let $K_n\Subset \Omega^o$ and $K_n\uparrow \Omega^o$. Then, let $\chi_n\in C_0^\infty(\Omega^o)$ and $\chi_n\equiv 1$ on $K_n$. Then, applying Proposition \ref{prop:prop}, to $\chi_n u$, we have that $\WF_h(\chi_n u)\cap K_n\times \re^d$ is invariant under the leaves generated by $H_{\Im p}$ and $H_{\Re p}.$ But, this is true for all $n$, so, letting $n\to \infty$, we obtain the result.
\end{proof}

\section{A Carleman Type Estimate}
\label{sec:Carleman}

We now prove a Carleman type estimate for $(P,\Omega)$. This will be used in the following sections to restrict the essential support of quasimodes.

Observe that for $\varphi\in C^\infty(\overline{\Omega})$, we have 
\begin{equation}\label{eqn:Pphi}P_{z,\varphi}:=e^{\varphi/h}P_ze^{-\varphi/h}=\sum(hD_{x_j}+i\partial_{x_j}\varphi)^2-\la X,\partial \varphi\ra +i\la X, hD\ra-z\end{equation}
with Weyl symbol 
\begin{equation}\label{eqn:Pphisymb} 
p_{z,\varphi}(x,\xi)=|\xi|^2-\la X+\partial \varphi,\partial \varphi\ra +i\la X+2\partial \varphi,\xi\ra -z.\end{equation}

Then, $P_{z,\varphi}=A+iB$ where $A$ and $B$ are formally self adjoint and have 
$$A=(hD)^2-\la X+ \partial \varphi,\partial \varphi\ra -\Re z,\quad B=\la X,hD\ra+\sum_j(\partial_{x_j}\varphi\composed hD_{x_j}+hD_{x_j}\composed \partial_{x_j}\varphi )-\Im z$$
with Weyl symbols
$$a=|\xi|^2-\la X+\partial \varphi,\partial\varphi\ra -\Re z,\quad b=\la X+2\partial \varphi,\xi\ra -\Im z.$$

%\noindent Also, note that 
%$$[A,B]=hi^{-1}\left(4\la \partial^2\varphi hD,hD\ra +\la \partial^2\varphi(X+2\partial\varphi),X+2\partial\varphi\ra +4hi^{-1}\la \partial \Delta \varphi,hD\ra +h^2i^{-2}\Delta \Delta\varphi\right).$$

Next, let $u\in C^\infty(\overline{\Omega})$ with $u|_{\partial \Omega}=0$, $Pu=v$, $u_1:=e^{\varphi/h}u$, and $v_1=e^{\varphi/h}v$. Then, we compute
\begin{eqnarray}
\|v_1\|^2&=&\left((A+iB)u_1,(A+iB)u_1\right)\nn
&=&\|Au_1\|^2+\|Bu_1\|^2+i\left[\left(Bu_1,Au_1\right)-\left(Au_1,Bu_1\right)\right]\nonumber 
\end{eqnarray}
Now, observe that, since $B$ is a first order differential operator that is formally self adjoint, and $u|_{\partial\Omega}=0$, 
\begin{equation}\label{eqn:firstOrder}
\left(Au_1,Bu_1\right)=\left(BAu_1,u_1\right).
\end{equation}
Next,
\begin{equation}\label{eqn:boundaryTerms}
\left(Bu_1,Au_1\right)=\left(ABu_1,u_1\right)-h^2\left(Bu_1,\partial_\nu u_1 \right)_{\partial\Omega}.\end{equation}
But, on $\partial \Omega$
$$B=\frac{h}{i}\la 2\partial \varphi+X,\nu\ra\partial_\nu +B'$$
where $B'$ acts along $\partial \Omega$. Hence, 
$$\left(Bu_1,\partial_\nu u_1\right)_{\partial\Omega}=\frac{h}{i}\left(\la 2\partial \varphi+X,\nu\ra \partial_\nu u_1,\partial_\nu u_1\right)_{\partial \Omega}$$
and we have 
\begin{equation}\label{eqn:prelimCarl}
\|v_1\|^2=\|Au_1\|^2+\|Bu_1\|^2+i\left([A,B]u_1,u_1\right)-h^3\left(\la 2\partial \varphi+X,\nu\ra \partial_\nu u_1,\partial \nu u_1\right)_{\partial \Omega}.\end{equation}

Next, we compute 
$$\{a,b\}=4\la \partial^2\varphi \xi,\xi\ra+\la \partial^2\varphi (X+2\partial\varphi),X+2\partial\varphi\ra.$$
Thus, choosing $\varphi=\e\psi$ with $\partial^2\psi$ positive definite, we have 
\begin{equation}
\label{eqn:poissonBracket}\{a,b\}\geq C\e|\xi|^2+C\e|X+2\e\partial\psi|^2\geq C\e(|\xi|^2+|X|^2)+O(\e^2).\end{equation}
Now,  $i[A,B]=h\{a,b\}^w+\e h^2r^w$, where %$r=-4\la\partial\Delta\psi,\xi\ra+hi^{-1}\Delta\Delta\psi$ 
$r\in S^1$. Hence, for $\delta >0$ small enough and independent of $h$, $h$ small enough, and $0<\e<\delta$ (here $\e$ may depend on $h$),  we have 
\begin{equation}
\label{eqn:commutator}
i[A,B]=Ch\e(-h^2(\partial^2\psi)^{ij}\partial^2_{x_ix_j}+ f(x))+\e O_{H_h^1\to L^2}(h^2)
\end{equation}
where $ f\geq C> 0$ and $\partial^2\psi>C\geq 0$. Hence, by an integration by parts,
$$i\left([A,B]u_1,u_1\right)\geq Ch\e(\|hDu_1\|^2+\|u_1\|^2).$$
Combining this with (\ref{eqn:prelimCarl}), noting that, on $\partial\Omega$, $\partial_\nu u_1= e^{\frac{\e\psi}{h}}\partial_\nu u$, and, letting  $\Gamma_+$ and $\partial\Omega_-$ be as in (\ref{eqn:OmegaBound}),
 we have,
\begin{multline*}
-h^3\left(\la 2\e\partial\psi+ X,\nu\ra  e^{\frac{\e\psi}{h}}\partial \nu u,e^{\frac{\e\psi}{h}}\partial \nu u\right)_{\partial\Omega_-}+Ch\e\|e^{\frac{\e\psi}{h}}u\|^2_{H_h^1}\\
\leq \|e^{\frac{\e\psi}{h}}P_zu\|^2+h^3\left(\la 2\e\partial\psi+ X, \nu\ra e^{\frac{\e\psi}{h}}\partial \nu u, e^{\frac{\e\psi}{h}}\partial \nu u\right)_{\Gamma_+}
\end{multline*}

Now, note that a similar proof goes through if 
\begin{equation}
\label{eqn:alternatePsi}\psi(x)=\psi_1(\la x,X\ra)
\end{equation}
 where $\psi_1$ has  $\partial^2\psi_1>0$. In this case \eqref{eqn:poissonBracket} reads 
  $$\{a,b\}\geq C\e|\la X,\xi\ra|^2+C\e|X+2\e\partial\psi|^2\geq C\e(|\la X,\xi\ra|^2+|X|^2)+O(\e^2)$$
 and \eqref{eqn:commutator} reads
 $$
 i[A,B]=Ch\e(-h^2(\partial^2\psi)^{ij}\partial^2_{x_ix_j}+ f(x))+\e O_{H_X\to L^2}(h^2)
 $$
 where $\|u \|_{H_X}=\|\la X,hD\ra u\|_{L^2}+\|u\|_{L^2}.$ After this observation, we obtain the following lemma,

\begin{lemma}\label{lem:Carleman}
Let $u\in C^\infty(\overline{\Omega})$ $u|_{\partial \Omega}=0$, $\psi\in C^\infty(\overline{\Omega})$ either have 
\begin{enumerate}
\item $\psi$ is locally strictly convex ($\partial^2\psi$ is positive definite), or
\item $\psi$ is as in \eqref{eqn:alternatePsi}.
\end{enumerate}
 Then, there exists $\delta>0$ independent of $h$ small enough such that for $0<\e\leq \delta$ ($\e$ possibly depending on $h$), and $0<h<h_0$, we have
\begin{multline}
\label{eqn:CarlemanEstimate}
-h^3\left(\la 2\e\partial\psi+ X,\nu\ra e^{\frac{\e\psi}{h}}\partial \nu u, e^{\frac{\e\psi}{h}}\partial \nu u\right)_{\partial\Omega_-}+Ch\e\|e^{\frac{\e\psi}{h}}u\|^2_{W}\\\leq \|e^{\frac{\e\psi}{h}}P_zu\|_{L^2}^2+h^3\left(\la 2\e\partial\psi+ X, \nu\ra e^{\frac{\e\psi}{h}}\partial \nu u, e^{\frac{\e\psi}{h}}\partial \nu u\right)_{\Gamma_+}
\end{multline}
where if $\psi$ satisfies
\begin{enumerate}
\item $\|\cdot\|_W=\|\cdot\|_{H_h^1}$
\item $\|\cdot\|_W=\|\cdot\|_{H_X}=\|\cdot \|_{L^2}+\|\la X,hD\ra \cdot \|_{L^2}.$
\end{enumerate} 
\end{lemma}

\noindent Lemma \ref{lem:Carleman} easily extends to $u\in H_h^2$ with $u|_{\partial \Omega}=0$. 

\section{Essential Support of Quasimodes}
\label{sec:QuasimodesOnIlluminated}
In this section, we prove part 2 of Theorem \ref{thm:PseudoSpec}.
\subsection{No quasimodes on the boundary of the pseudospectrum}
\label{sec:boundPseud}
Let $z_0\in \partial \Lambda(P,\Omega).$ We use a small weight to conjugate $P$ as in (\ref{eqn:Pphi}) such that $p_\varphi$ is elliptic. For simplicity, we again assume $X=e_1$ and hence $\Re z_0=(\Im z_0)^2$. Using (\ref{eqn:Pphisymb}), let $\e>0$ and $\partial \varphi=-\e X$ (i.e. $\varphi=-\e\la X,x\ra +C$). Then, using the fact that $X=e_1$, we have
$$p_{z,\varphi}(x,\xi)=|\xi|^2+ (1-\e)\e+i(1-2\e)\xi_1-z_0.$$

Then, $p_{z,\varphi}=0$ implies that 
$$\xi_1=\frac{\Im z_0}{1-2\e}\quad |\xi '|^2+\recip{(1-2\e)^2}\left((\Im z_0 )^2-\Re z_0-\Re z_0(-4\e +4\e^2)\right)+\e(1-\e)=0.$$
But,
\begin{multline*}|\xi '|^2+\recip{(1-2\e)^2}\left((\Im z_0 )^2-\Re z_0-\Re z_0(-4\e +4\e^2)\right)+\e(1-\e)\\\geq \recip{(1-2\e)^2}(4\e\Re z_0(1-\e))+(1-\e)\e>0\end{multline*}
for $\e$ small enough. We now show that $$|p_{z,\varphi}|\geq c\e\la \xi\ra^2$$
for $\e$ small enough. The fact that $|p_{z,\varphi}|\geq c\la \xi\ra^2$ for $|\xi|>>1$, is clear. Thus, we only need to check that $|p_{z,\varphi}|>c\e$. Let $\xi_1=(\Im z_0+\gamma)/(1-2\e).$ Then, choose $\gamma =\delta\e$ 
\begin{eqnarray}
|p_{z,\varphi}(x,\xi)|&\geq&|\xi '|^2+\recip{(1-2\e)^2}(\gamma^2+2\gamma \Im z_0+4\e \Re z_0-4\e^2\Re z_0)+\e-\e^2\nn 
&\geq&(1+O(\e))(\delta^2\e^2+2\delta\e\Im z_0+4\e \Re z_0-4\e^2\Re z_0)+\e +O(\e^2)\nn
&\geq&\e(2\delta\Im z_0+4\Re z_0+1)-O(\e^2)\geq c'\e\nonumber 
\end{eqnarray}
for $\delta$ small enough independent of $\e$ and for $\e$ small enough. 

Therefore, by Proposition \ref{prop:Calderon}, if $u|_{\partial \Omega}=0$, we have that 
$$\e\|e^{\frac{\e \psi}{h}}u\|_{H_h^2}\leq  C\|e^{\frac{\e\psi}{h}}P_z u\|_{L^2}.$$
Thus, if $u$ is a quasimode for $z_0$, choosing $\e=h\log h^{-1}$, 
$$\|u\|_{L^2}\leq C(h^N\log h^{-1})^{-1}O(h^\infty)=O(h^\infty),$$
a contradiction. Hence, there are no quasimodes for $z_0\in \partial \Lambda(P,\Omega)$. 

Thus, we have proved 
\begin{lemma}
\label{lem:NoPseudBoundaryQuasi}
Suppose $\Re z_0=|X|^{-2}(\Im z_0)^2$. Then there are no quasimodes of $(P,\Omega)$ for $z_0$. 
\end{lemma}

\noindent {\bf Remark:} This argument can be adjusted slightly to give that if $d(z_0, \partial \Lambda(\Omega,P))=O(h),$ then there are no quasimodes for $z_0\in \partial \Lambda(P,\Omega).$

\subsection{No Quasimodes Away from the Illuminated Boundary}
To finish the proof of Theorem \ref{thm:PseudoSpec}, we will need the following elementary lemma. (The proof follows \cite[Section 6.3.2]{E}.) Let $Q(x,hD)u=-h^2\partial_j(c^{ij}\partial_iu) +\la a(x), h\partial u\ra +b(x)u$. 
\begin{lemma}\label{lem:boundaryReg}
Suppose that $\partial \Omega\in C^1$ and that $a,b,c^{ij}\in C^\infty(\re^d;\complex)$ with $c^{ij}\xi_i\xi_j\geq C|\xi|^2$ uniformly in $\Omega$. Then there exists $C>0$ such that for all $u\in C^\infty(\Omega)\cap H_0^1(\Omega)$ we have 
$$\|u\|_{H_h^2(\Omega)}\leq C(\|Qu\|_{L^2(\Omega)}+\|u\|_{L^2(\Omega)}).$$
\end{lemma}
\begin{proof}
Using a partition of unity and change of coordinates, we can assume without loss of generality that $\Omega=B(0,1)\cap \{x_1\geq 0\}$. Then, let $\chi\in C^\infty(\overline{\Omega})$ 
with $\chi\equiv 1$ on $V:=\{|x|<\recip{2}\}$ and $\chi\equiv 0$ on $|x|>\frac{3}{4}$. 
 Next, let $v=-h^2\partial_k\chi^2\partial_k\bar{u}$ for $k=2,...,d$. Then, $v\in H_h^1$ with $v|_{\partial \Omega}=0$ and hence
$$\int  c^{ij}h\partial_i u h\partial_j v =\int Qu v-\la a(x), h\partial u\ra v-b(x)uv.$$
Now,
\begin{align*}
\int  c^{ij}h\partial_i u h\partial_j v&=\int -c^{ij} h\partial_i u h^2\partial_j\partial_k(\chi^2h\partial_k\bar{u})\\
&=\int \left[c^{ij}h^2\partial_i\partial_ku+h\partial_k(c^{ij})(h\partial_iu)\right]\left[\chi^2h^2\partial_j\partial_k\bar{u}+2\chi (h\partial_j\chi)(h\partial_k\bar{u})\right]\\
&\geq \int \chi^2(1-C\e)|h^2\partial \partial_k u|^2-h^2\e^{-1}|h\partial_ku|^2\\
&\geq \int \recip{2}\chi^2|h^2\partial \partial_ku|^2-Ch^2|h\partial_ku|^2.
\end{align*}
Then, 
\begin{align*}
\left|\int Qu v-\la a(x), h\partial u\ra v-b(x)uv\right|&\leq \int (|Qu|+|\la a(x), h\partial u\ra|+|b(x)u|)|v|\\
&= \int (|Qu|+|\la a(x), h\partial u\ra|+|b(x)u|)\left|2\chi h\partial_k\chi h\partial_k\bar{u}+\chi^2h^2\partial_k^2\bar{u}\right|\\
&\leq \int C|Qu|^2+C|h\partial u|^2+C|u|^2+Ch^2|hDu|^2+C\e |\chi^2h^2\partial_k^2u|^2\nonumber
\end{align*}
Thus, 
$$\|h^2\partial_k\partial u\|_{L^2(V)}\leq C(\|Qu\|_{L^2(\Omega)}+\|u\|_{H_h^1(\Omega)})$$
for $k=2,...,d$. Now, for $k=1$, we note that 
$$h^2\partial_1^2u=(c^{11})^{-1}\left(Qu+\sum_{(i,j)\neq (1,1)}c^{ij}h^2\partial_i\partial_ju-\sum_{i,j=1}^{d}(\partial_jc^{ij})\partial_iu-\la a(x),h\partial u\ra-b(x)u\right)$$
where $(c^{11})^{-1}$ is well defined by the positive definiteness of $c^{ij}.$
Thus, 
$$|h^2\partial_1^2u|\leq C\left(\sum_{i=2}^{d}|h^2\partial_i^2u|+|hD u|+|u|+|Qu|\right)$$
and we have 
$$\|u\|_{H_h^2(V)}\leq C(\|Qu\|_{L^2(\Omega)}+\|u\|_{H_h^1(\Omega)})$$
and the result follows from \cite[Theorem 7.1]{EZB} and its proof.
\end{proof}
We apply the above lemma to obtain the following,

\begin{lemma}
\label{lem:EStoSupp}
Suppose that $u$ has $u|_{\partial \Omega}=0$, $\|u\|_{L^2}=1$, and $\ES_h(u)\cup \ES_h(P_zu)\subset A$ and $P_zu=O_{L^2}(1)$. Then, for any $U$ with $A\Subset U$, and $\chi\in C^\infty(\overline{\Omega})$ with $\chi\equiv 1$ on $U$, 
$$\|(1-\chi)u\|_{H_h^2}=O(h^\infty).$$
In particular, if $u$ is a quasimode for \eqref{eqn:generic} with $\ES_h(u)\subset A$, then, for any $U$ with $A\Subset U$, there is a quasimode $u_1$ with supp $u_1\subset U$. 
\end{lemma}
\begin{proof}
Let $A\Subset U_0\Subset U\Subset U_1\Subset B$. Let $\chi\in C^\infty(\overline{\Omega})$ have $\chi\equiv 1 $ on $U$ and supp $\chi \Subset U_1$.  Let $\chi_0=\chi$ and for $i=1,...$ let $\chi_i\in C^\infty(\overline{\Omega})$ have supp $\chi_i \subset B\setminus U_0$ and have $\chi_i\equiv 1$ on supp $\partial \chi_{i-1}$. Then, by Lemma \ref{lem:boundaryReg}
\begin{align*}
\|(1-\chi)u\|_{H_h^2}&\leq C(\|P_z(1-\chi)u\|_{L^2}+\|(1-\chi)u\|_{L^2})\\
&\leq C\|(1-\chi)P_zu\|_{L^2}+C\|[P_z,\chi]u\|_{L^2}+O(h^\infty)\\
&=O(h^\infty)+\|[P_z,\chi]u\|_{L^2} \leq O(h^\infty)+Ch\|\chi_1 u\|_{H_h^1}.
\end{align*}

But, using the same argument again, we have that since $\chi_n\equiv 0$ on $U_0$ for all $n$, 
$$\|\chi_{n-1}u\|_{H_h^1}\leq \|(1-(1-\chi_{n-1}))u\|_{H_h^2}\leq O(h^\infty)+Ch\|\chi_{n} u\|_{H_h^1}$$
Hence, by induction, for all $N>0$,
$$\|(1-\chi)u\|_{H_h^2}\leq O(h^\infty)+C_Nh^N\|\chi_Nu\|_{H_h^1}.$$
But, by Lemma \ref{lem:boundaryReg}, since $Pu=O_{L^2}(1)$, $u=O_{H_h^2}(1)$ and hence
$$\|(1-\chi)u\|_{H_h^2}\leq O(h^\infty)+C_Nh^N\|u\|_{H_h^1}=O(h^\infty)$$
as desired.

To prove the second claim observe that if $u$ is a quasimode,
$$\|P_z\chi u\|_{L^2}\leq \|P_z(1-\chi)u\|_{L^2}+\|P_zu\|_{L^2}=O(h^\infty)$$
since $(1-\chi)u=O_{H_h^2}(h^\infty).$
\end{proof}

We now apply the above  lemma to restrict the essential support of quasimodes.
\begin{lemma}
\label{lem:isectIllum}
If $u$ is a quasimode for \eqref{eqn:generic} then $\ES_h(u)\cap\overline{ \partial\Omega_+}\neq \emptyset.$
\end{lemma}
\begin{proof}
Suppose that $u$ is a quasimode for \eqref{eqn:generic} and $u$ has 
$$\ES_h(u)\cap \overline{\partial\Omega_+}=\emptyset.$$
Then, by Lemma \ref{lem:EStoSupp}, we may assume that $u$ is supported away from $\overline{\partial \Omega_+}$.

Now, applying Lemma \ref{lem:Carleman} with $\e=h\log h^{-1}$ and $\psi=\la X,x\ra^2$, we have 
\begin{multline*}-h^3\left(\la 2\e\partial\psi+ X,\nu\ra  e^{\frac{\e\psi}{h}}\partial \nu u,e^{\frac{\e\psi}{h}}\partial \nu u\right)_{\partial\Omega_-}+Ch\e\|e^{\frac{\e \psi}{h}}u\|^2_{L^2}\\\leq O(h^\infty)+h^3\left(\la 2\e\partial\psi+ X,\nu\ra e^{\frac{\e\psi}{h}}\partial \nu u, e^{\frac{\e\psi}{h}}\partial \nu u\right)_{\partial\Omega_0}.\end{multline*}
But, since $\partial\psi=2\la X,x\ra X$, the term on $\partial\Omega_0$ vanishes and, hence, we have 
$$Ch\e\|e^{\frac{\e\psi}{h}}u\|^2_{L^2}= O(h^\infty).$$
Hence, $u=O_{L^2}(h^\infty)$ and there are no quasimodes concentrating away from $\partial\Omega_+$ - i.e.
$$\ES_h(u)\cap\overline{\partial\Omega_+}\neq \emptyset.$$ 
\end{proof}

\subsection{Characterization of the Essential Support of Quasimodes}

In order to use Lemma \ref{lem:Carleman} to characterize $\ES_h(u)$ for quasimodes, we would like to construct a set $A$ with $\Gamma_+\subset A$ and a weight function $\psi$ such that for any $U\subset \overline{\Omega}$ separated from $A$, there exists $\e>0$ such that $\sup_A\psi<\inf_U\psi-\e$. Since $\psi$ must be locally convex in $\overline{\Omega}$ to apply Lemma \ref{lem:Carleman}, any set $A$ with this property must be relatively convex inside $\overline{\Omega}$ (Recall that relative convexity is defined in Definition \ref{def:ConvSub}.). 

\subsubsection{Preliminaries on Relatively Convex Sets}
Let $B$ be a bounded set and $A$ be convex relative to $B$. We wish to determine whether there is a smooth locally strictly convex function (inside $B$) with $\partial A$ as a level set. 

\begin{lemma}
\label{lem:gauge}
Let $A$ be a closed and relatively convex set inside $B$, a bounded set. Then there is a function $g_A$ that is locally convex inside $B$ and has $g_A|_{A}\equiv 0$, $g_A(x)>0$ for $x\notin A$. 
\end{lemma}
\begin{proof}
First, define the epigraph of a function $f$ as follows. 
$$\text{epi}(f)=\{(x,\mu)\in B\times \re:\mu\geq f(x)\}.$$
We show that a function $f$ is locally convex in $B$ if and only if its epigraph is relatively convex in $B\times \re$. 

Suppose that $f$ is locally convex in $B$. Then, for every $x,y\in B$ with $L_{x,y}\subset B$, $f(tx+(1-t)y)\leq tf(x)+(1-t)f(y).$ (Here $L_{x,y}$ is as in \eqref{eqn:lineSeg}) Therefore, if $(x,\mu),(y,\nu)\in$ epi$(f)$, then $(tx+(1-t)y,t\mu+(1-t)\nu)\in \text{ epi}(f).$

Now, suppose that $\text{ epi}(f)$ is relatively convex in $B\times \re$. Then, suppose that $x,y\in B$ with $L_{x,y}\subset B$. Then, let $f(x)=\mu$ and $f(y)=\nu$. Then, $t(x,\mu)+(1-t)(y,\nu)\in \text{ epi}(f).$ Hence,
$$f(tx+(1-t)y)\leq t\mu+(1-t)\nu=tf(x)+(1-t)f(y)$$
and $f$ is locally convex in $B$. 

Now, we determine the epigraph of the function $g_A$. First, let 
$$G=A\times[0,\infty)\cup B\setminus A\times[1,\infty).$$
Then, let $G_A=\ch_{B\times \re}(G)$. Observe that since $A$ is relatively convex in $B$, $A\times [0,\infty)$ is relatively convex in $B\times \re$. Now, by Carath\'{e}odory's Theorem, any point in $\ch(G)$ can be written as the convex combination of at most $d+2$ points in $G$. Since $A\times [0,\infty)$ is relatively convex in $B\times \re$, and $\ch_{B\times \re}(G)\subset \ch(G)$, any point in $G_A\setminus(A\times[0,\infty))$ is representable as a convex combination of $d+2$ points, at least one of which is in $B\setminus A\times[1,\infty)$. 

Suppose $x\notin A$ and $(x,\nu)\in G_A$. Then, $d(x,A)>0$ since $A$ is closed and 
$$(x,\nu)=\sum_{i=1}^{d+2}t_i(x_i,\nu_i)$$
where, for some $r>0$, $x_1,...,x_r\notin A$. Hence, $\nu_1,...,\nu_r\geq 1$. Relabel $(x_i,\nu_i)$ $i=1,...,r$ so that $t_1=\max(t_1,...,t_r)$. Then, since $B$ is bounded there exists $R>0$ such that $B\subset B(0,R)$ and hence $t_1> \tfrac{d(x,A)}{(d+2)R}$. Therefore 
$$\nu\geq \frac{d(x,A)}{(d+2)R}>0 .$$ 
 Thus letting 
$$g_A(x)=\inf\{y:(x,y)\in G_A\},$$ 
$g_A$ is locally convex in $B$ with $g_A>0$ on $B\setminus A$ and $g_A=0$ on $A$. 
\end{proof}

\begin{corol}
\label{cor:convex}
Let $B$ and $B_1\Supset B$ be bounded sets. Let $A\subset B_1$ be closed and convex relative to $B_1$. Then there exists $\psi\in C^\infty(\overline{B})$ strictly locally convex in $B$  such that for all $W$ with $A\Subset W$, $\sup_A \psi<\inf_{B\setminus W}\psi$.
\end{corol}
\begin{proof}
Let $d(B_1,B)=2\delta$. Then, let $\varphi_\e\in C_0^\infty(\re^d)$ be a nonnegative approximate identity family with support contained in $B(0,\delta)$ and define $f_A^\e=g_A*\varphi_\e$ where $g_A$ was constructed in Lemma \ref{lem:gauge}. (Here we extend $g_A$ off of $B_1$ by 0.) Then, $f_A^\e\to g_A$ uniformly on bounded sets. Also, $f_A^\e$ is smooth. To see that $f_A^\e$ is locally convex inside $B$, observe that for $x,y\in B$ with $L_{x,y}\subset B$,
\begin{multline*}f_A^\e(tx+(1-t)y)=\int \varphi_\e(z)g_A(tx+(1-t)y-z)dz\\\leq \int \varphi_\e(z)(tg_A(x-z)+(1-t)g_A(y-z))dz=tf_A^\e(x)+(1-t)f_A^\e(y)\end{multline*}
by the local convexity of $g_A$ inside $B_1$ and the nonnegativity of $\varphi_\e$. Finally, to make a locally strictly convex approximation of $g_A$, define $g_A^\e:=f_A^\e+\e|x|^2$. Then, $g_A^\e\to g_A$ uniformly on bounded sets, and $g_A^\e\in C^\infty(\overline{B})$ with $g_A^\e$ locally strictly convex inside $B$. 
\end{proof}

\noindent{\bf Remark: }Although we have not constructed a smooth locally convex function with level set $\partial A$, we have one that has a level set which is uniformly arbitrarily close.

We also need a few more properties of relatively convex sets 
\begin{lemma}
\label{lem:closureIsConvex}
Suppose that $A\subset B$ is relatively convex in $B$, $B$ open and bounded. Then, $\overline{A}$ is relatively convex in $B$.
\end{lemma}
\begin{proof}
Let $x, y\in \overline{A}$ such that $L_{x,y}\subset B$. Then, there are sequences $x_n\to x$ and $y_n\to y$ with $x_n,y_n\subset A$. We need to show that $L_{x,y}\subset\overline{A}.$ For $0\leq \lambda \leq 1$, We have that 
$$|\lambda(x_n-x)+(1-\lambda)(y_n-y)|\leq |x_n-x|+|y_n-y|.$$
But, since $L_{x,y}$ is compact and $B$ is open, there is $\e>0$ such that 
$$\{z:d(z,L_{x,y})<\e\}\subset B$$
and hence, we have that for $n$ large enough $L_{x_n,y_n}\subset B$. But, since $A$ is relatively convex, this implies $L_{x_n,y_n}\subset A$ and hence for $0\leq \lambda \leq 1$
$$\lim_{n\to \infty} \lambda x_n+(1-\lambda)y_n=\lambda x+(1-\lambda )y\in \overline{A}.$$
\end{proof}
\begin{lemma}
\label{lem:ConvIsecSupset}
We have that 
$$\bigcap_{B\Subset B_1}\ch_{B_1}(A)=\ch_{\overline{B}}(A).$$
%and if $B$ is compact,
%$$\bigcap_{A\Subset A_1}\ch_{B}(A_1)=\ch_{B}(\overline{A}).$$
\end{lemma}
\begin{proof}
Let 
$$\mcal{C}=\{C:A\subset C\,, x,y\in C\,, L_{x,y}\Subset B_1\text{ for all }B_1\Supset B \text{ implies } L_{x,y}\subset C\}.$$
Then,
$$\bigcap_{B\Subset B_1}\ch_{B_1}(A)=\bigcap_{C\in \mcal{C}}C.$$
But, if $L_{x,y}\nsubseteq \overline{B}$, then $L_{x,y}\nsubseteq B_1$ for some $B_1\Supset B$. Hence, 
$$\mcal{C}=\{C:A\subset C\,, x,y\in C\,, L_{x,y}\Subset \overline{B} \text{ implies } L_{x,y}\subset C\}$$
and the result follows.

%{\color{red}Finish this proof}
%To see the second result, observe that since $B$ is compact, if $L_{x,y}\nsubseteq B$, then 
%$$\sup_{z\in L_{x,y}}d(x,B)>0.$$
%Hence, there exist neighborhoods of $x$ and $y$, $U$ and $V$ such that for all $z\in U$, $w\in V$, $L_{z,w}\nsubseteq B$. Hence, since $\overline{A}$ is compact, for all $x,y\in \overline{A}$, either $L_{x,y}\subset A$ or 
\end{proof}

\subsubsection{Application to Quasimodes}
\label{sec:convexRestriction}
We now apply the above results on relatively convex sets to quasimodes. 
\begin{lemma}
\label{lem:convexRestriction}
Let $u\in L^2$ have $u|_{\partial \Omega}=0$, $\ES_h(u)\subset A$, $P_zu=O_{L^2}(1)$, and $\ES_h(P_zu)\subset B$. Then, for all $A_1\Supset A$ and $B_1\Supset B$,
$$\ES_h(u)\subset \overline{\ch_{\overline{\Omega}}((A_1\cap\Gamma_+)\cup B_1)}.$$ 
In particular, if $u$ is a quasimode for \eqref{eqn:generic}, then for all $A_1\Supset A$
$$\ES_h(u)\subset \overline{\ch_{\overline{\Omega}}(A_1\cap \Gamma_+)}.$$
\end{lemma}
\begin{proof}
Choose $\Omega_1$ open with $\overline{\Omega}\Subset \Omega_1$, $A_1$ closed with $A\Subset A_1$, and $B_1$ closed with $B\Subset B_1$. Let $F=\overline{\ch_{\Omega_1}((A_1\cap\Gamma_+)\cup B_1)}.$ Then, by Lemma \ref{lem:closureIsConvex}, $F$ is relatively convex in $\Omega_1$.  Let $U\subset\overline{\Omega}$ such that $d(U,F)>0$. By Corollary \ref{cor:convex}, there exists $\psi\in C^\infty(\overline{\Omega})$ strictly convex in $\overline{\Omega}$ such that for some $\delta>0$, $\sup_F\psi<\inf_U\psi-\delta$. 

\noindent {\bf Remark:} The regions of $\psi>0$ and $\psi=0$ are shown in Figure \ref{f:Omega}.

We have, by Lemma \ref{lem:Carleman} that, for $\e_0$ small enough independent of $h$, and $h$ small enough, for $\e(h)<\e_0$,

$$Ch\e\|e^{\frac{\e\psi}{h}}u\|^2_{H_h^1}\leq C\|e^{\frac{\e\psi}{h}}P_zu\|^2+Ch^3\left(\la 2\e\partial\psi+ X, \nu\ra e^{\frac{\e\psi}{h}}\partial \nu u, e^{\frac{\e\psi}{h}}\partial \nu u\right)_{\Gamma_+}.$$

Now, suppose that $u$ has $\ES_h(u)\subset A$, $\ES_h(P_z u)\subset B$, and let $\e=\gamma h\log h^{-1}$. Then for $A\Subset A_1$ and $B\Subset B_1$ Lemma \ref{lem:EStoSupp} gives that
up to $O_{H_h^2}(h^\infty)$ supp $u\subset A_1\cup B_1$. Thus, 
$$h^3\left(\la 2\e\partial\psi+ X, \nu\ra e^{\frac{\e\psi}{h}}\partial \nu u, e^{\frac{\e\psi}{h}}\partial \nu u\right)_{ \Gamma_+}\leq O(h^\infty)+h^3\left(\la 2\e\partial\psi+ X, \nu\ra e^{\frac{\e\psi}{h}}\partial \nu u, e^{\frac{\e\psi}{h}}\partial \nu u\right)_{(A_1\cup B_1)\cap \Gamma_+}.$$ Hence, we have 
$$Ch\e\|e^{\frac{\e\psi}{h}}u\|^2_{L^2}\leq O(h^\infty)+\|e^{\e\psi/h}P_zu\|_{L^2(B_1)}+h^3\left(\la 2\e\partial\psi+ X, \nu\ra e^{\frac{\e\psi}{h}}\partial \nu u, e^{\frac{\e\psi}{h}}\partial \nu u\right)_{(A_1\cup B_1)\cap \Gamma_+}.$$
But, by Lemma \ref{lem:boundaryReg}, $\|u\|_{H_h^2}=O(1)$. Hence, 
$$h^3(\partial_\nu u, \partial_\nu u)_{\partial\Omega}\leq C h^3\|u\|_{H^{3/2}}^2=C\|u\|_{H_h^{3/2}}^2\leq C\|u\|^2_{H_h^2}=O(1).$$
Thus, we have that 
$$Ch\e\inf_U e^{\frac{2\e\psi}{h}}\|u\|^2_{L^2(U)}\leq Ch\e\|e^{\frac{\e\psi}{h}}u\|^2_{L^2(U)}\leq Ch\e\|e^{\frac{\e\psi}{h}}u\|^2_{L^2}\leq O(h^\infty)+C\sup_{(A_1\cap\Gamma_+)\cup B_1}e^{\frac{2\e\psi}{h}}.$$
But, $\inf_U\psi\geq \delta+\sup_F\psi$ and we have 
$$C_\gamma h^{2-2\gamma \delta}\log h^{-1}\|u\|_{L^2(U)} = Ch\e e^{\frac{2\e\delta}{h}}\|u\|^2_{L^2(U)}\leq O(h^\infty)+C.$$
Hence, letting $\gamma\to \infty$, we have that $\|u\|_{L^2(U)}=O(h^\infty)$ as desired. 

Thus, $u$ cannot have essential support away from $F$. That is for any $A_1\Supset A$ and $B_1\Supset B$,
\begin{equation}
\nonumber
\ES_h(u)\subset \bigcap_{\Omega\Subset \Omega_1}\overline{\ch_{\Omega_1}((A_1\cap \Gamma_+)\cup B_1)}=\overline{\ch_{\overline{\Omega}}((A_1\cap \Gamma_+)\cup B_1)}%=\overline{\ch_{\overline{\Omega}}((\overline{A}\cap\Gamma_+)\cup \overline{B})}=\ch_{\overline{\Omega}}((\overline{A}\cap\Gamma_+)\cup \overline{B}).
\end{equation}
Here, equality of the two sets follows from Lemma \ref{lem:ConvIsecSupset}.
 The second claim follows from the fact that a quasimode has $\ES_h(Pu)=\emptyset$. 
\end{proof}
\noindent {\bf Remark:} Observe that if $\Gamma_+\subset A$, then the second part of Lemma \ref{lem:convexRestriction} gives that for quasimodes $$\ES_h(u)\subset \overline{\ch_{\overline{\Omega}}(\Gamma_+)}.$$
\subsection{Characterization of the Interior Wavefront set of a Quasimode}
We wish to determine the possible essential support of a quasimode. To do this we first need the following simple lemmas
\begin{lemma}
\label{lem:relateESWF}
For a solution to (\ref{eqn:generic}), 
$$\pi_x(\WF_h(u))\cap\Omega^o =\ES_h(u)\cap\Omega^o.$$
\end{lemma} 
\begin{proof}
Say $x_0\in \pi_x(\WF_h(u))\cap \Omega^o$. Then it is clear that $x_0\in \ES_h(u).$ 

Now, suppose $x_0\notin\pi_x(\WF_h(u))\cap\Omega^o$.  Let $K>0$ such that $|p_z(x,\xi)|\geq C\la \xi\ra^2$ for $|\xi|\geq K$ and $x\in \Omega^o$. Let $U$ be a neighborhood of $x_0$ such that for all $\chi\in C_0^\infty(U\times\{|\xi|\leq 2K\})$,  $$\|\chi^wu\|_{L^2}=O(h^\infty).$$
Such a neighborhood, $U$ exists by the compactness of $\{|\xi|\leq 2K\}$ and \cite[Theorem 8.13]{EZB}. 

Let $x_0\in V\Subset U$, $\varphi\in C_0^\infty(U)$ with $\varphi\equiv 1$ on $V$, and $\psi\in C^\infty(\re^d)$ with supp $\psi\subset\{|\xi|\leq 2K\}.$

\noindent To complete the proof, we need only show that there is a $V$ such that 
$$\|(1-\psi(\xi))^w\varphi u\|_{L^2(V)}=O(h^\infty).$$

\noindent To see this, let $\psi_1\in C_0^\infty(\re^{d})$ have $\psi\equiv 1 $ on supp $\psi_1$ and supp $\psi_1 \subset\{|\xi|\leq 2K\}$, let $\varphi_1\in C_0^\infty(U)$ with $\varphi_1\equiv 1$ on supp $\varphi$, and finally let $\varphi_2\in C_0^\infty(U)$ with $\varphi_2\equiv 1$ on supp $\varphi_1$. 

Then, observe that 
$$(1-\psi_1)^2|p_z|^2\varphi_2^2\varphi_1^2\geq \gamma\la \xi\ra^4$$
on supp $(1-\psi)\varphi$. 
Hence, by the Sharp G$\mathring{\text{a}}$rding inequality, 
$$\|\varphi_2P_z\la hD\ra^{-2}\la hD\ra^2(1-\psi)^w\varphi u\|_{L^2}^2\geq \gamma^2\|(1-\psi)^w\varphi u\|_{L^2}^2-Ch\|\la hD\ra^2(1-\psi)^w\varphi u\|_{L^2}^2.$$
But, by Lemma \ref{lem:boundaryReg} 
$$\|(1-\psi)^wu\|_{H_h^2}\leq C(\|(1-\psi)^w\varphi u\|_{L^2}+\|P_z(1-\psi)^w\varphi u\|_{L^2})$$
and we have that 
$$\frac{\gamma}{2}\|(1-\psi)^w\varphi u\|_{L^2}\leq C\|\varphi_2P_z(1-\psi)^w\varphi  u\|_{L^2}$$
But, 
$$\varphi_2P_z(1-\psi)^w\varphi u=\varphi_2P_z u+\varphi_2P_z((1-\psi^w)\varphi-1) u=O_{L^2}(h^\infty).$$ since $P_zu=O_{L^2}(h^\infty)$ and $\varphi_2P_z(1-(1-\psi^w)\varphi)=c^w$ with supp $c\subset U\times\{|\xi|\leq 2K\}$.
\end{proof}

\begin{figure}[htbp]
\includegraphics[width=5in]{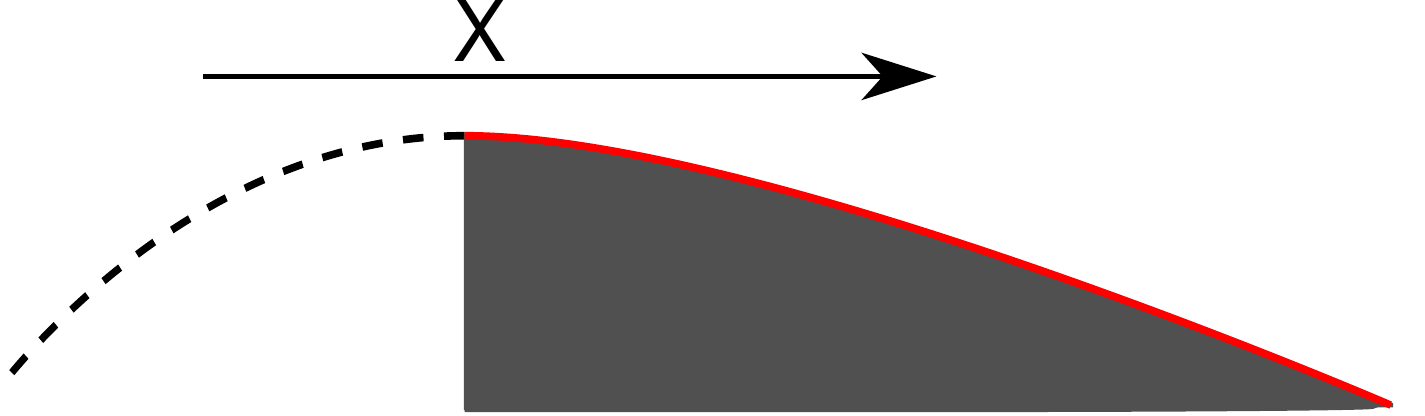}
\begin{center}
\caption{\label{f:outsideConvex}
The figure shows a piece of $A\cap \Omega$ where $A$ is a plane tangent to $X$. A portion of $\partial\Omega_-$ is shown in the black dashed line and a portion of $\partial\Omega_+$ is shown in the red line. The shaded region represents the convex hull of $\Gamma_+$. The non-shaded region is a portion of $\Omega \cap A$ that is not contained in $\ch(\Gamma_+)$.}
\end{center}
\end{figure}

\begin{lemma}
\label{lem:outsideConvex}
Let $\partial \Omega\in C^1$. Then for any plane $A$ with $X$ tangent to $A$ and $A\cap\Omega^o\neq \emptyset$, we have 
\begin{equation}
\nonumber
\left(\Omega^o\setminus \ch(\Gamma_+)\right)\cap A\neq \emptyset.
\end{equation}
\end{lemma}
\begin{proof}
For simplicity, we again assume $X=e_1$. By the compactness of $\Gamma_+$ we can choose $x\in \Gamma_+\cap A$ such that $\pi_1(x)\leq \pi_1(y)$ for all $y\in \Gamma_+\cap A$ where $\pi_1$ is projection onto the first component. Then,
$$\pi_1\left(\ch(\Gamma_+)\cap A\right)\subset [\pi_1(x),\infty).$$
We show that there is a $z\in \Omega^o\cap A$ with $\pi_1(z')<\pi_1(x)$ and hence that 
$z'\notin \ch(\Gamma_+)\cap A.$

\noindent {\bf Remark:} The regions of interest are shown in Figure \ref{f:outsideConvex}.

Suppose $x\in \partial \Omega_+$. Then, $\la e_1,\nu(x)\ra>0$ and hence there is $z\in \Omega^o\cap A$ with $\pi_1(z')<\pi_1(x)$. Now, suppose that $x\in \partial \Omega_0$. Then, $e_1$ is tangent to $\partial\Omega\cap A$ at $x$. Hence, there is a $z\in \partial\Omega\cap A$ with $\pi_1(z)<\pi_1(x)$. But, this implies that there is a $z'\in \Omega^o\cap A$ with $\pi_1(z')<\pi_1(x).$
\end{proof}

We now finish the proof of part (2) of Theorem \ref{thm:PseudoSpec}.
\begin{proof}

Let $u$ be a quasimode. Observe that if $(p-z)(x_0,\xi_0)=0$ and $\xi_0\neq (\Im z,0,...,0)$, then $\Re z>(\Im z)^2$. Hence, by Lemma \ref{lem:WFSubsetCharSet}, and Corollary \ref{cor:Propagation}, if 
$$x_0\in \Omega^o,\quad (x_0,\xi_0)\in \WF_h(u),$$
then, there exists a plane $A$ tangent to $e_1$ with $x_0\in A$ such that
$$\pi_x(\WF_h(u))\supset A\cap \Omega^o.$$
But, $\ES_h(u)$ is closed and $\ES_h(u)\cap \Omega^o=\pi_x(\WF_h(u)\cap (\Omega^o\times\re^d))$, hence
$$\ES_h(u)\supset A\cap \overline{\Omega}.$$
Together with Lemma \ref{lem:convexRestriction}, this gives for $\Gamma_+\Subset A_1$
$$A\cap \overline{\Omega}\subset \ES_h(u)\subset \overline{\ch_{\overline{\Omega}}(A_1)}.$$
But, notice that
$$ \overline{\ch_{\overline{\Omega}}(A_1)}\subset \overline{\ch(A_1)}.$$
Hence, since $A_1\Supset \Gamma_+$ was arbitrary, 
$$A\cap \overline{\Omega}\subset \ES_h(u)\subset \bigcap_{A_1\Supset \Gamma_+}\overline{\ch(A_1)}=\ch(\Gamma_+)$$
since $\Gamma_+$ is compact.
Therefore we have a contradiction of Lemma \ref{lem:outsideConvex}.
Putting this together with Lemma \ref{lem:WFSubsetCharSet}, we have
$$\WF_h(u)\cap(\Omega^o\times \re^d)\subset (\Omega^o\times \{(\Im z,0,0,...,0)\})\cap p_z^{-1}(0).$$

Now, note that if $\xi=(\Im z, 0,...,0)$, and $p(x,\xi)-z=0$, then $\Re z=(\Im z)^2$ and hence $z\in \overline{\Lambda(P,\Omega)}.$ 
Thus, except for $z\in \partial \overline{\Lambda(P,\Omega)}$, 
$$\WF_h(u)\cap(\Omega^o\times\re^d)=\emptyset.$$
But, we have shown in Lemma \ref{lem:NoPseudBoundaryQuasi} that there are no quasimodes for $z\in \partial \overline{\Lambda(P,\Omega)}$. Hence quasimodes cannot have wave front set in the interior of $\Omega$.

So, using Lemma \ref{lem:relateESWF}, we have 
\begin{equation}\label{eqn:interiorES}\ES_h(u)\cap \Omega^o=\pi_x(\WF_h(u)\cap(\Omega^o\times \re^d))=\emptyset.\end{equation}
Thus, $\ES_h(u)\subset \partial \Omega$. 

To finish the proof of Theorem \ref{thm:PseudoSpec}, we apply Lemma \ref{lem:convexRestriction} with $\ES_h(u)\subset \partial\Omega$ to obtain
\begin{equation}\label{eqn:convexQuasiES}\ES_h(u)\subset \overline{\ch_{\overline{\Omega}}(\Gamma_+)}.\end{equation}
\noindent Putting \eqref{eqn:interiorES} and \eqref{eqn:convexQuasiES} together, we have that quasimodes cannot concentrate away from the intersection of the  $\overline{\Omega}$ convex hull of the glancing and illuminated boundary with the boundary - i.e. 
 $$\ES_h(u)\subset \partial\Omega\cap\overline{\ch_{\overline{\Omega}}(\Gamma_+)}$$
as desired. 
\end{proof}

\subsection{Further Localization}
We now apply Lemma \ref{lem:Carleman} locally to obtain further information about the essential support of quasimodes -- we prove parts (4) and (5) of Theorem \ref{thm:PseudoSpec}. 

We will need the following lemma.
\begin{lemma}
\label{lem:localizeQuasi}
Let $\chi\in C^\infty(\overline{\Omega})$, then for any quasimode $u$ of \eqref{eqn:generic}, any $A,B$ with $\supp\partial\chi\cap \partial\Omega\Subset B$ and $\supp \chi \cap \partial\Omega\Subset A$, we have
$$\ES_h(\chi u)\subset \overline{\ch_{\supp \chi\cap \overline{\Omega}}((\Gamma_+\cap A)\cup B)}\cap \partial\Omega.$$
\end{lemma}
\begin{proof}
Let $u$ be a quasimode for \eqref{eqn:generic}, $\chi\in C^\infty(\overline{\Omega})$. 
Now, let $W$ have supp $\partial \chi \Subset W$ and let $U_1$ be a neighborhood of $W\cap \partial \Omega$.
 Then, let $\chi_1\in C^\infty(\overline{\Omega})$ with $\chi_1\equiv 1$ on $U_1$. Then,
$$\|(1-\chi_1)P_z\chi u\|_{L^2}=O(h^\infty)+\|(1-\chi_1)[P_z,\chi]u\|_{L^2}.$$
Now, by Lemma \ref{lem:EStoSupp}, and the fact that $\ES_h(u)\subset \partial\Omega$,
$$\|(1-\chi_1)[P_z,\chi]u\|_{L^2}\leq C\|(1-\chi_1)u\|_{H_h^1(W)}=O(h^\infty).$$
 Hence, since $U_1$ was an arbitrary neighborhood of supp $\partial\chi \cap\partial \Omega$,  $$\ES_h(P_z\chi u)\subset \text{supp } \partial \chi \cap \partial \Omega.$$

Then, observe that $\chi u$ is a function on $\Omega_1=\text{ supp }\chi \cap \overline{\Omega}$ with 
$$\chi u|_{\partial\Omega_1}=0\,,\quad \ES_h(\chi u)\subset \partial\Omega\cap \text{ supp }\chi\,,\quad \ES_h(P_z\chi u)\subset \text{ supp } \partial \chi \cap \partial \Omega.$$

\noindent Hence, applying Lemma \ref{lem:convexRestriction} to $\chi u$ on $\Omega_1$, and using the fact that $\ES_h(u)\subset \partial \Omega$, we have for every $B\Supset \text{ supp } \partial \chi\cap\partial\Omega $ and every $A\Supset \text{ supp } \chi\cap \partial\Omega$, 
\begin{equation}
\label{eqn:cutoff}\ES_h(\chi u)\subset \overline{\ch_{\text{supp }\chi\cap\overline{\Omega}}((\Gamma_+\cap A)\cup B)}\cap \partial \Omega.\end{equation}
\end{proof}

We now use Lemma \ref{lem:localizeQuasi} to finish the proof of Theorem \ref{thm:PseudoSpec}.

\begin{figure}[htbp]
\includegraphics[width=3.5in]{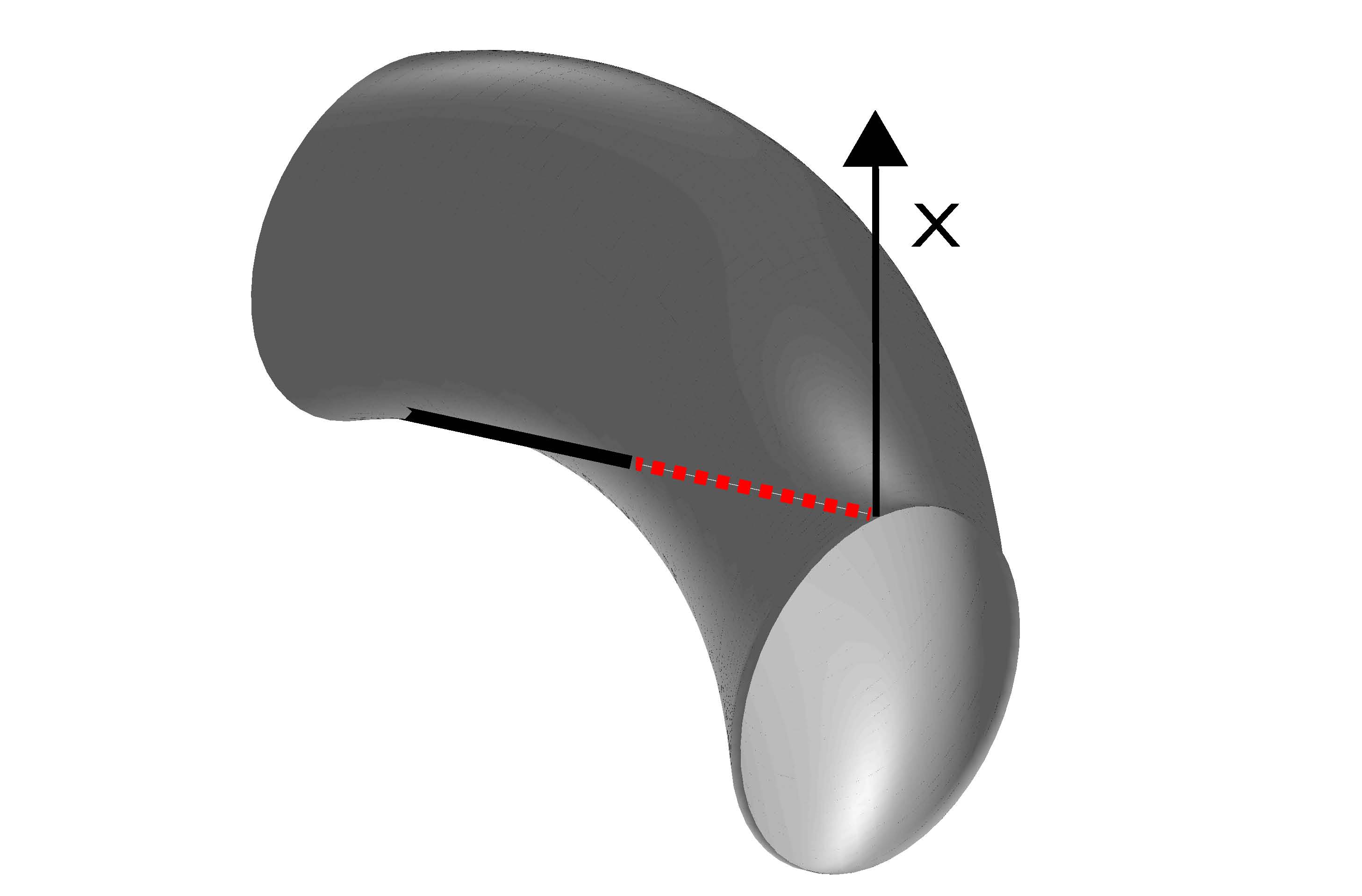}
\begin{center}
\caption{\label{f:badShape}
The figure shows a piece of a domain $\Omega\subset \re^3$. A portion of $\partial\Omega_-$ is shown in the black line and a portion of $\partial\Omega_+$ is shown in the dashed red line. Notice that for any point in the portion of $\partial\Omega_-$ shown, $\partial\Omega_+$ can be reached along a straight line lying entirely inside the boundary.  This example shows that in dimensions higher than 2 we cannot hope to make an argument similar to that used to prove the last part of Theorem \ref{thm:PseudoSpec}.}
\end{center}
\end{figure}

\begin{proof}
For simplicity, assume $X=e_1$. 
To prove the first part of the proposition, suppose that $\partial\Omega$ is either strictly concave or strictly convex at $x_0\in \partial\Omega_-$. Then there exists $\chi\in C^\infty(\overline{\Omega})$ such that $\chi\equiv 1$ in a neighborhood of $x_0$, $\supp \chi \cap \partial\Omega \Subset \partial\Omega_-$ and for $x\in \supp \partial\chi\cap \partial\Omega$, 
$$|\pi_1(x)-\pi_1(x_0)|>\delta.$$

Then there exists $A$ with $\partial\Omega_-\supset A\Supset \supp \chi \cap \partial\Omega$ and $B\Supset \supp \partial\chi \cap \partial\Omega$ such that for $x\in B$
$$|\pi_1(x)-\pi_1(x_0)|>\delta/2.$$
Hence, 
$$x\in \overline{\ch_{\supp \chi\cap \overline{\Omega}}((\Gamma_+\cap A)\cup B)}\cap \partial\Omega=\overline{\ch_{\supp \chi\cap \overline{\Omega}}(B)}\cap \partial\Omega$$
implies 
$$|\pi_1(x)-\pi_1(x_0)|>\delta/2$$
and by Lemma \ref{lem:localizeQuasi}, $x_0\notin \ES_h(\chi u)$. Thus, since $\chi\equiv 1 $ in a neighborhood of $x_0$, $x_0\notin \ES_h(u).$

Now, suppose $\Omega\subset \re^2$. Then choose $x_0\in \partial\Omega_-$ and let $\gamma:[-1,1]\to \partial \Omega$ be a curve defining $\partial\Omega$ with $\gamma(0)=x_0$, $\gamma(-1),\gamma(1)\in \partial\Omega_0$, and $\gamma((-1,1))\subset \partial\Omega_-$. Defining
$$t_{\pm}:=\inf\{t:\gamma '(\pm t)\neq \gamma '(0)\},$$
we have $|t_{\pm}|<1$ since if not, then $\la X,\gamma ' (\pm 1)\ra\neq 0$. Then, there exists $\e>0$, such that for $r>0$ small enough
$$x_0\notin B(\gamma(t_--\e),r)\cup B(\gamma(t_++\e),r)=:W_1\cup W_2,\quad (W_1\cup W_2)\cap \Gamma_+= \emptyset,$$
and there exists $\delta >0$ such that 
\begin{equation}
\label{eqn:awayBoundary}\inf_{(z_1,z_2)\in W_1\times W_2}\sup_{s\in [0,1]} d(sz_1+(1-s)z_2,\partial\Omega)>\delta.
\end{equation}

Let $\chi\in C^\infty(\overline{\Omega})$ have $\chi\equiv 1$ in a neighborhood of $x_0$,
$$\supp \partial\chi \cap \partial\Omega\Subset W_1\cup W_2, \quad \supp \chi \cap \partial\Omega \Subset\partial\Omega_-,$$
and
$$\supp \chi \subset \{x\in \overline{\Omega}:d(x,\partial\Omega)<\delta/2\}.$$
\noindent Then, letting $u$ be a quasimode for \eqref{eqn:generic}, and applying Lemma \ref{lem:localizeQuasi}
$$\ES_h(\chi u)\subset \overline{\ch_{\supp \chi\cap \overline{\Omega}}(W_1\cup W_2)}=\overline{W_1\cup W_2}.$$
Here, the last equality follows from \eqref{eqn:awayBoundary} and the convexity of $B(x,r)$. 
Hence, $x_0\notin \ES_h(\chi u)$ and we have $x_0\notin \ES_h(u).$ But, $x_0\in \partial \Omega_-$ was arbitrary. Therefore, $\ES_h(u)\subset \Gamma_+$ as desired. 
\end{proof}
\noindent{\bf Remark:} Figure \ref{f:badShape} shows an example of why we cannot make a similar argument in dimensions larger than 2.

\section{Instability in an Evolution Problem}
\label{sec:Instabilty}
Our approach to obtaining blow-up of (\ref{eqn:evolutionEq}) will follow that used by Sandstede and Scheel in \cite{Sand03} and that by the author in \cite{Ga}. We first demonstrate that, from small initial data, we obtain a solution that is $\geq 1$ on a translated ball in time $t_1=O(1)$. We then use the fact that the solution is $\geq 1$ on this region to demonstrate that, after an additional $t_2=O(h)$, the solution to the equation blows up. 

First, we prove that there exists initial data so that the solution to (\ref{eqn:evolutionEq}) is $\geq 1$ in time $O(1)$. Let $\varphi_t:=\exp(-tX)$ denote the flow of $i\la X,D\ra$. Note that for the purposes of Theorem \ref{thm:evolution}, we do not need to assume that $X$ is constant. 
\begin{lemma}
\label{lem:O(1)ndim}
Fix $\mu>0$, $\alpha<\mu$, $0<\e\leq\recip{2}(\mu-\alpha)$, and $(x_0,a,\delta)\in \re^d\times \re^+\times \re^+$ such that both $\varphi_t(B\left(x_0,2a\right))\subset \Omega$ for $0\leq t\leq 2 \delta$ and $\varphi_t$ is defined on $B\left(x_0,3a\right)$ for $0\leq t<2\delta$. Then, for each 
$$0<h<h_0$$
 where $h_0$ is small enough, there exists
$$u_0(x)\geq 0,\quad \|u_0\|_{C^k}\leq \exp\left(-\recip{C_kh}\right),\quad k=0,1,...$$ 
and $0<t_1<\delta$ so that the solution to (\ref{eqn:evolutionEq}) with initial data $u_0$ satisfies $u(x,t)\geq 1$ on $x\in \varphi_{t}(B(x_0,a))$ for $t_1\leq t<\delta$.
\end{lemma}

\noindent\begin{proof}
The proof of this lemma follows that in \cite[Lemma 3]{Ga} except we no longer need to control the size of the potential. Instead, we show that the ansatz satisfies Dirichlet boundary conditions on $\partial\Omega$. 

Let $\upsilon$ solve 
\begin{equation}
\label{eqn:linProb}
(h\partial_t+P(x,hD)-\mu)\upsilon=0,\quad \upsilon(x,0)=\upsilon_0,\quad \upsilon|_{\partial\Omega}=0.
\end{equation}

Let $w_0:\re^d\to \re$ and define $O:=\{x:w_0>0\}$. We make the following assumptions on $w_0$,
\begin{equation}\label{eqn:wAssume1}
w_0\geq 0,\quad\|w_0\|_{C^k}\leq \exp(-\recip{C_kh}),\quad w_0\in C(\re^d)
\end{equation}
\begin{equation}\label{eqn:wAssume2}
w_0\in C^\infty(\overline{O}),\quad \text{supp }w_0\subset B(x_0,2a),\quad
w_0>\exp\left(-\frac{\delta}{2h}\right)\text{ on }B(x_0,a),
\end{equation}
\begin{equation}\label{eqn:wAssume3}
\partial O\text{ is smooth},\quad -\Delta w_0 (x)\leq Cw_0(x)-\beta\text{ for }x\in O\text{ and }0<h<h_0.
\end{equation}
where $C^\infty(\overline{O})$ are smoothly extendible functions on $O$. We refer the reader to \cite[Lemma 3]{Ga} for the construction of such a function.

Define $w:[0,2\delta)\times \re^d\to \re$ by 
 $$w:=\begin{cases}
	\exp\left(\frac\alpha h t\right)w_0(\varphi_{t}(x))&\text{ where }\varphi_t\text{ is defined},\\
0&\text{else}.\end{cases}$$
Since supp $w\subset B(x_0,2a)$ and $\varphi_{t}$ is defined on $B(x_0,2a)\times [0,2\delta)$, $w$ is continuous. We proceed by showing that $w$ is a viscosity subsolution of $(\ref{eqn:linProb})$ in the sense of Crandall, Ishii, and Lions \cite{Lions}. 

First, we show that $w$ is a subsolution on $O_t:=\varphi_t(O)$ for $t<\delta$. 
\begin{align*}
hw_t+P(x,hD)w-\mu w&=hw_t-h^2\Delta w+ih\la X ,D\ra w-\mu w \\
&=(\alpha -\mu )w-h^2\Delta w\\ 
&\leq  \exp\left(\frac \alpha ht\right)
\left((\alpha - \mu )w_0\right)-h^2\Delta w,
\end{align*}
(Here, we evaluate all instances of $w_0$ at $\varphi_t(x)$.)
Now, by Taylor's formula, for $x\in \Omega$, $\varphi_t(x)=x+O(t)$ (with similar estimates on $x$ derivatives). Hence $-\Delta\left[w_0(\varphi_t(x))\right]=-\Delta w_0(\varphi_t(x))+O(t)$.  We have $t<\delta$, and $-\Delta w_0\leq Cw_0-\beta$ on $O$. Therefore, for $\delta$ small enough, $-\Delta w\leq Cw_0 $. Hence, for $h$ small enough independent of $0<\delta <\delta_0$,
$$hw_t+P(x,hD)w-\mu w\leq \exp\left(\frac{\alpha}{h}t\right)\left( \alpha-\mu+Ch^2\right)w_0\leq 0$$
Now, since for $t<\delta$, $\text{supp } w\subset \Omega$ we have that $w$ is a subsolution on $O_t$ for $t<\delta$ and $h$ small enough. Next, observe that on $(\re^d\setminus \overline{O_t})$, $w\equiv 0$ and hence is a subsolution of (\ref{eqn:linProb}) on this set as well.  

Finally, we need to show that $w$ is a subsolution on $\partial O_t:= \varphi_t(\partial O)$. We refer the reader to the proof of \cite[Lemma 3]{Ga} for this. Lastly, observe that since $\varphi_t(B(x_0,2a))\subset \Omega$ for $t<2\delta$, we have that for $t<\delta$, $w|_{\partial \Omega}=0$. Together with the previous arguments, this shows that $w$ is a viscosity subsolution for (\ref{eqn:linProb}) on $t<\delta$. 

Now, by an adaptation of the maximum principle found in \cite[Section 3]{Lions} to parabolic equations, any solution, $\upsilon$ to (\ref{eqn:linProb}) with initial data $\upsilon_0>w_0$ has $\upsilon\geq w$ for $t<\delta$. Now, suppose $u_1$ solves 
$$h\partial_tu_1+(P-\mu)u_1=|u_1|^p,\,u_1|_{t=0}=\upsilon_0\geq 0.$$
Then, $u_1$ is a supersolution for \eqref{eqn:linProb} and hence has $u_1\geq\upsilon\geq 0$. But this implies that in fact $u:=u_1$ solves \eqref{eqn:evolutionEq} with initial data $\upsilon_0$. Therefore, $u\geq\upsilon\geq w$ for $t<\delta$ and hence, since for $t>\frac{\delta}{2}$, $w(x,t)\geq 1$ on $\varphi_t(B(x_0,a))$, we have the result.
\end{proof}

\noindent{\bf Remark:} To obtain a growing subsolution it was critical that $\mu>0$. This corresponds precisely with the movement of the pseudospectrum of $(-(P-\mu),\Omega)$ into the right half plane.

Now, we demonstrate finite time blow-up using the fact that in time $O(1)$ the solution to (\ref{eqn:evolutionEq}) is $\geq 1$ on an open region. Again, the proof of Theorem \ref{thm:evolution} follows that in \cite[Theorem 1]{Ga} except we replace the need to control the size of the potential with the requirement that the solution be 0 on $\partial \Omega$. 

\begin{proof}
Let $u_0(x)$ and $t_1$ be the initial data and time found in Lemma \ref{lem:O(1)ndim} with $(a, x_0,\delta)$ such that $\varphi_t$ is defined  on $B(x_0,a)$, $\varphi_{t}\left(B\left(x_0,a\right)\right)\subset \Omega$ for $t\in [0,\delta]$, and $t_1<\delta$. Then, $u(x,t_1)\geq 1$ on $\varphi_t\left(B\left(x_0,a\right)\right)$. 

Now, let $\Phi\in C_0^\infty(\re)$ be a smooth bump function with $\Phi(y)=1$ on $|y|\leq 1$, $0\leq \Phi \leq 1$, supp $\Phi\subset (-2,2)$, and $\Phi '',\Phi'\leq C\Phi^{1/p}$. Define $\chi:\re^d\to \re$ by $\chi(y):=\Phi\left(2a^{-1}|y|\right).$ To see that such a function $\Phi$ exists, observe that when $\Phi>c>0$, the inequality can easily be arranged by adjusting $C$. Then, notice that the function $e^{-1/x}$ has 
$$\left(e^{-1/x}\right)''=e^{-1/x}(x^{-4}-2x^{-3})\leq Ce^{-1/(px)}$$
and
$$\left(e^{-1/x}\right)'=e^{-1/x}x^{-2}\leq Ce^{-1/(px)}$$
for $x$ small enough. 

Next, let $y'=\varphi_{t}(x_0+y)$ and let
$$v(y,t):=\chi(y)u(y',t).$$
Then, we have that 
$$hv_t=h^2\Delta v+\mu v+v^p-2h^2\la\nabla\chi ,\nabla u\ra -h^2u\Delta \chi +(\chi -\chi ^p)u^p.$$
Finally, define the operations, $[f]$ and $( f,g)$ by
$$[f]:=\negint_{B(0,a)}f(y)dy\quad ( f,g ):=\negint_{B(0,a)}\la f(y),g(y)\ra dy.$$
(Here, $\int\negthickspace\negthickspace\negthickspace-$ denotes averaging.)

Then,  
\begin{eqnarray}
h[v]_t&=&h^2[\Delta v]+\mu [v]+[v^p]-h^2\left(\Delta \chi ,u\right) -2h^2\left[\la \nabla\chi,\nabla u\ra\right] +\left( \chi -\chi ^p,u^p\right)  \nn 
&\geq &\mu [v]+[v^p]+h^2\left( \Delta \chi ,u\right) +\left( \chi -\chi ^p,u^p\right) \label{eqn:intPartsnd}
\end{eqnarray}
Here, (\ref{eqn:intPartsnd}) follows from integration by parts, and the fact that $\nabla \chi=0$ at $|y|= a$.

\noindent We will later need that $[v^p]\geq [v]^p$. To see this use H\"{o}lder's inequality as follows 
$$[v]^p=\recip{|B(0,a)|}\int_{B(0,a)}v^p\leq \int_{B(0,a)}\frac{v^p}{|B(0,a)|}\left(\int_{B(0,a)}\recip{|B(0,a)|}\right)^{p-1}=[v^p].$$
We will also need an estimate on $\left( \Delta\chi , u\right)$. Following \cite[Section 4]{Ga}, we obtain
\begin{align}\left|\left( \Delta\chi ,u\right)\right|&=\left|\frac{1}{C_da^d}\int_0^a\int_{S^{d-1}}\frac{2(d-1)}{ar}\left[\Phi'(2a^{-1}r)+4a^{-2}\Phi''(2a^{-1}r)\right]r^{d-1}u(r\phi)dS(\phi)dr\right|\nonumber \\
&\leq C\left|\int_0^a\int_{S^{d-1}}\Phi^{1/p}r^{d-1}udS(\phi)dr\right|\nonumber\\
&\leq C\int \chi^{1/p}u\leq C\negint(1+\chi u^p)\leq C'+C\negint\chi u^p\label{eqn:cubicnd}
\end{align}
where $C'$ and $C$ do not depend on $h$.

Now, we have 
\begin{eqnarray}
h[v]_t&\geq&\mu [v]+[v^p]+h^2\left( \Delta \chi ,u\right) +\left( \chi -\chi ^p,u^p\right)\nn
&\geq& \mu [v]+[v^p]-O(h^2)+\left( (1- O(h^2))\chi -\chi^p,u^p\right)\nn 
&\geq &\mu [v]+[v^p]-O(h^2)-O(h^2)[v^p] \label{eqn:holder}\\
&\geq&\mu [v] +(1-O(h^2))[v]^p-O(h^2)\nonumber
\end{eqnarray}
Here, (\ref{eqn:holder}) follows from the fact that $\chi \leq 1$ and $[v^p]\geq [v]^p$. Note that these equations are satisfied for $t<\delta$ since $\varphi_t(B(x_0,a))\subset \Omega$ for $t<\delta$.

We have that $[v](t_1)\geq 1/4$ and $\mu>0$. Then, by Lemma \ref{lem:O(1)ndim}, there exists $\gamma>0$ independent of $h$ such that for $h$ small enough and $t_1\leq t\leq t_1+\gamma$,
$$h[v]_t\geq \frac{\mu}{2}[v]+\recip{2}[v]^p.$$ 
But, the solution to this equation with initial data $[v](0)\geq 1/4$ blows up in time $t_2=O(h)$. Hence, so long as $t_1+t_2<\min(\delta,t_1+\gamma)$ and $h$ is small enough, $[v]$ blows up in time $t_1+t_2$. Observe that since $t_1<\delta$, $0\leq t_1+t_2=t_1+O(h)<\min(\delta,t_1+\gamma)$ for $h$ small enough. Thus, the solution to \eqref{eqn:evolutionEq} blows up in time $\delta$. 
\end{proof}

\section{Application to Hitting Times for Diffusion Processes}
\label{sec:meanHit}
Let $\Omega\subset \re^d$ be a bounded domain with $C^\infty$ boundary. Then, define the stochastic process 
\begin{equation}
\label{eqn:diffusionProcess}dX_t=b(X_t)+\sqrt{2h}dB_t
\end{equation}
\begin{figure}[htbp]
\includegraphics[width=3.5in]{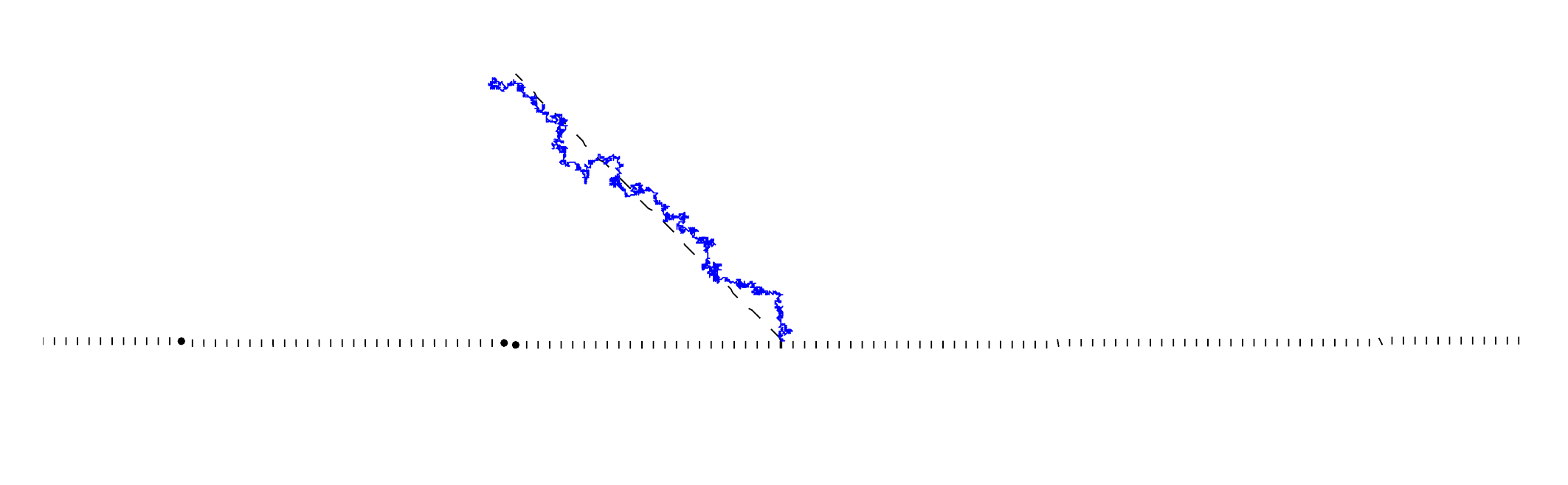}
\begin{center}
\caption{\label{f:diffuse}
The figure shows a sample of the diffusion process $X_t$ with $b(X_t)=(-\recip{\sqrt{2}},\recip{\sqrt{2}})$, $h=10^{-4}$, and initial condition $X_0=(0,h)$. The dotted line shows the boundary of a disk tangent to $y=0$ of radius $1/2$ and the dashed line shows the path of the ode with no noise. The boundary is shown for $y<10h$.}
\end{center}
\end{figure}
where $B_t$ is Brownian motion and $b\in C^\infty(\re^d;\re^d)$. (Figure \ref{f:diffuse} shows an example path for $X_t$.) Let $Y_t=X_{ht}.$ Then, $Y_{t}$ solves
$$dY_t=hb(Y_t)+\sqrt{2}hdB_{t}\quad Y_0=x_0$$
and it is a standard result of probability theory \cite[Section 1.5]{FreidlinWentzell} that the operator 
\begin{equation}\label{eqn:diffPDE}
L:=-(hD)^2+i\la b,hD\ra
\end{equation}
is associated to $Y_t$ in the sense that if
\begin{equation}
\label{eqn:diffusion}\begin{cases}
(-L-\lambda)u=((hD)^2+i\la -b,hD\ra -\lambda)u=f&\text{ in }\Omega,\\
u|_{\partial \Omega}=0,\end{cases}\end{equation}
then $u$ has
\begin{equation}u(x)=\mathbb{E}_x\int_0^{\tau_Y} f(Y_t)e^{\lambda t}dt.\label{eqn:meanHittingEqn}\end{equation}
were $\mathbb{E}_{x}$ denotes the expected value given that $X_0=x$.

Next, define the first hitting times, by 
\begin{equation}
\label{eqn:exitTimes}
\tau_Y:=\inf\{t\geq 0: Y_t\in \partial \Omega\}\quad \tau_X:=\inf\{t\geq 0:X_t\in \partial \Omega\}=h\tau_Y.
\end{equation}

\noindent Let $\lambda_1(-L)$ denote the principal eigenvalue of $-L$. We prove the following proposition,
\begin{prop}
\label{prop:exitTime}
Let $\Omega\subset\re^d$ be bounded with $\partial \Omega\in C^\infty$. Let $X_t$ and $\tau_X$ be defined as in \eqref{eqn:diffusionProcess} and \eqref{eqn:exitTimes} respectively. Then, for each $x_0\in \partial \Omega_+$ (where $X=-b\in C^\infty(\re^d;\re^d)$ in \eqref{eqn:OmegaBound} -- that is $\la b,\nu\ra <0$.), and
$$\begin{cases}0<\lambda <\lambda_1(-L)\, , \lambda\neq \la b(x_0),\nu(x_0)\ra^2/4&d\geq 2\\
0<\lambda <\min(\lambda_1(-L),\la b(x_0),\nu(x_0)\ra^2/4)&d=1
\end{cases},\quad \lambda \text{ independent of }h.$$
There exists $\gamma>0$ such that for all $N\geq 1$ and $x(h)$ with $$\gamma h>|x(h)|>c_1h^N,\quad \la x(h),\nu(x_0)\ra <-c_2<0,\quad x_0+x(h)\in \Omega\,,\text{ for }0<h<h_0$$
 there exists $C>0$ such that for $h$ small enough,
$$h\log \mathbb{E}_{x_0+x(h)}e^{\lambda\tau_X/h}\geq Ch\log h^{-1}.$$
Moreover, if $\partial \Omega$ and $b$ are real analytic near $x_0$, there is a $\delta>0$ such that 
$$h\log \mathbb{E}_{x_0+x(h)}e^{\lambda\tau_X/h}\geq \delta$$
and such that for every $\alpha>1$ and $\e>0$, there exists $c_\alpha>0$ ($c_\alpha$ depending only on $\alpha$) and a function $s(h)>\delta-h^{1-\e}$ with  
\begin{equation}
\label{eqn:probBound}
\min(c_\alpha e^{-\alpha(s(h)-\delta)/h},1) \leq P\left(\tau_X\geq \frac{s(h)}{\lambda}\right)\leq P\left(\tau_X\geq \frac{\delta-h^{1-\e}}{\lambda}\right). 
\end{equation}
\end{prop}
\vbox{
\noindent{\bf Remarks:}
\begin{enumerate}
\item Notice that if $s(h)\leq  \delta$, for $0<h<h_0$, then by \eqref{eqn:probBound} we have that 
$$P(\tau_X\geq \delta/(2\lambda))\geq\min( c_a,1)$$
and hence, for $h$ small enough, that the first hitting time is larger than $\delta/(2\lambda)$ with uniformly positive probability. On the other hand, if $s(h)\to \infty$ or remains bounded but is $>\delta$, \eqref{eqn:probBound} gives control of the decay rate of $P(\tau_X\geq \delta /(2\lambda))$ as $h\to 0$.
\item If $b=-\nabla f$ for $f\in C^\infty(\re^d)$ and $|b|> 0$ in $\overline{\Omega}$, then $0<c<\lambda_1(-L)$ uniformly in $h$. Hence in these cases, there exist $\lambda$ as required by Proposition \ref{prop:exitTime}.
\item In fact, the proof gives that for all $\e>0$, there exists $h$ small enough so that we can take
$$c_\alpha =\frac{1-\alpha}{2-\alpha}(1-\e).$$
\end{enumerate}
}

\begin{proof}
Let $L$ be as in \eqref{eqn:diffPDE}. Then, for $\lambda< \lambda_1(-L)$, we have that the solution, $u$ to \eqref{eqn:diffusion} has \eqref{eqn:meanHittingEqn}. 

Now, by Proposition \ref{prop:constructQuasi}, if $0<\lambda$, there are quasimodes for \eqref{eqn:diffusion} that are concentrated near $x_0$ for $x_0$ in the subset of the boundary illuminated by $-b$. Let $\la \mu,\nu(x_0)\ra <0$, $\e(h)\leq \gamma h$. We change coordinates so that $\nu(x_0)=e_1$ and observe that near the point $x_0$, these quasimodes have 
\begin{align*}
|u(\mu\e(h)+x_0)|&=\left|e^{-c|\mu '|^2\e(h)^2/h}\right|\left|a(\mu\e(h)+x_0)e^{ic_1\mu_1\e(h)/h+O(\mu_1\mu '\e(h)^2/h)}\right.\\
&\qquad{}\qquad{}\qquad{}\qquad\left.{}-b(\mu\e(h)+x_0)e^{ic_2\mu_1\e(h)/h+O(\mu_1\mu '\e(h)^2/h)}\right|\\ 
&=(1+O(\e(h)))\left|(e^{ic_1\mu_1\e(h)/h}-e^{ic_2\mu_1\e(h)/h})\right|+O(\e(h))\geq C_{\mu}\e(h)/h
\end{align*}
Therefore, for $\gamma>0$ small enough, every $\mu$ with $\la \mu,\nu(x_0)\ra<0$ and $\e(h)\leq \gamma h$, we have $|u(\mu\e(h)+x_0)|\geq C\e(h)/h$. Now, applying this in \eqref{eqn:meanHittingEqn}, we have 
\begin{align*}
C\frac{\e(h)}{h}\leq \left|\mathbb{E}_{x_0+\mu\e(h)}\int_0^{\tau_Y} f(Y_t)e^{\lambda t}dt\right|&=\left|\mathbb{E}_{x_0+\mu\e(h)}\int_0^{\tau_Y} (-L-\lambda)u(Y_t)e^{\lambda t}dt\right|\\&\leq \|(-L-\lambda)u\|_{L^\infty}\recip{|\lambda|}\mathbb{E}_{x_0+\mu \e(h)}\left(e^{\lambda\tau_Y}-1\right).
\end{align*}
If $\partial \Omega$ and $b$ are real analytic near $x_0$, we have $\|(-L-\lambda )u\|_{L^\infty}=O(e^{-\delta/h})$ which yields
$$\frac{\e(h)}{h}e^{\delta/h}\leq \mathbb{E}_{x_0+\mu\e(h)}e^{\lambda \tau_Y}=\mathbb{E}_{x_0+\mu\e(h)}e^{\lambda \tau_X/h}$$
and if $\partial \Omega$ or $b$ is only $C^\infty$ near $x_0$, $\|(-L-\lambda)u\|_{L^\infty}=O(h^\infty)$ and hence, for all $N>0$ there exists $c_N$ such that 
\begin{equation}\nonumber c_N\e(h)h^{-N}\leq \mathbb{E}_{x_0+\mu\e(h)}e^{\lambda \tau_Y}=\mathbb{E}_{x_0+\mu \e(h)}e^{\lambda \tau_X/h}.\end{equation}
Thus, if there exists $N>0$ such that $\e(h)>Ch^N$, we have, possibly with a different $\delta$,
\begin{equation}E_{x_0+\mu \e(h)}e^{\lambda \tau_X/h}\geq \begin{cases}e^{\delta/h} &\partial\Omega, b\text{ analytic near }x_0,\\
c_Nh^{-N}&\partial \Omega, b\text{ }C^\infty\text{ near }x_0.\end{cases}\label{eqn:diffuseSmooth}\end{equation}
This gives the first two statements in Proposition \ref{prop:exitTime}.

\noindent{\bf Remark:} Notice also, applying the standard small noise perturbation results that can be found, for example, in \cite[Theorem 2.3]{FreidlinWentzell} to a domain $\Omega_\delta\supset \Omega$ with $B(x_0,\delta)\subset \Omega_\delta$, and defining $\tau^\delta_X$ the corresponding hitting time, that we have for some $C>0$
$$\mathbb{E}_{x_0+\mu \e(h)}e^{\lambda \tau_X/h}\leq \mathbb{E}_{x_0+\mu \e(h)}e^{\lambda \tau^\delta_X/h}\leq e^{C/h}.$$

We now prove the second part of the proposition. Compute, using the fact that $\tau_X\geq 0$ and making the change of variables $s=h\log x$, 
$$E_{x_0+\mu\e(h)}e^{\lambda \tau_X/h}=\int_0^\infty P\left(e^{\lambda \tau_X/h}\geq x\right)dx=\recip{h}\int_0^\infty e^{s/h}P\left(\tau_X\geq s\lambda^{-1}\right)ds.$$
Hence, in the analytic case, 
$$\recip{h}\int_0^\infty e^{s/h}P\left(\tau_X\geq s\lambda^{-1}\right)ds\geq e^{\delta/h}$$
and we have that 
$$\recip{h}\int_0^\infty e^{(s-\delta)/h}P\left(\tau_X\geq s\lambda^{-1}\right)ds\geq 1.$$
Now, making the change of variables $t=(s-\delta)/h$.
$$\int_{-\delta/h}^\infty e^tP\left(\tau_X\geq  \lambda^{-1}(ht+\delta)\right)dt=o(1)+\int_{-h^{-\e}}^\infty e^tP\left(\tau_X\geq  \lambda^{-1}(ht+\delta)\right)dt .$$
 Thus, choosing $g(t)\in L^1(e^tdt)$ with $g(t)>0$, we have for all $\e_1>0$ and  $h$ small enough
\begin{align*}
1-\e_1&\leq \int_{-h^{-\e}}^\infty e^tP\left(\tau_X\geq \lambda^{-1}(ht+\delta)\right)dt\\
&\leq \left\|P\left(\tau_X\geq \lambda^{-1}(ht+\delta)\right)\left[g(t)\right]^{-1}1_{t>-h^{-\e}}\right\|_{L^\infty}\int_{-h^{-\e}}^\infty g(t)e^tdt 
\end{align*}
and hence
$$\frac{1-\e_1}{\|g\|_{L^1(e^tdt)}}\leq \left\|P\left(\tau_X\geq \lambda^{-1}(ht+\delta)\right)\left[g(t)\right]^{-1}1_{t>-h^{-\e}}\right\|_{L^\infty}.$$
That is, letting $ht+\delta=s$, for all $\gamma>0$, and $\e>0$, there exists $s(h)\geq \delta-h^{1-\e} $ such that 
$$\left[\frac{1-\e_1}{\|g\|_{L^1(e^tdt)}}-\gamma\right]g\left(\frac{s(h)-\delta}{h}\right)\leq  P\left(\tau_X\geq s(h)\lambda^{-1}\right)\leq P\left(\tau_X\geq \frac{\delta-h^{1-\e}}{\lambda}\right).$$
Fixing $\alpha>1$, letting $g(t)=\min(e^{-\alpha t},1)$, and letting $\gamma=\e_1/\|g\|_{L^1(e^tdt)}$ gives the last part of Proposition \ref{prop:exitTime}.
\end{proof}

\end{document}